\documentclass[11pt]{article}

\usepackage{amstext,amssymb,amsmath,amsbsy}
\usepackage{hyperref}
\usepackage{amscd}
\usepackage{amsfonts}
\usepackage{indentfirst}
\usepackage{verbatim}
\usepackage{amsmath}
\usepackage{amsthm}
\usepackage{enumerate}
\usepackage{graphicx}
\usepackage{color}
\usepackage[OT1]{fontenc}
\usepackage[latin1]{inputenc}
\usepackage[english]{babel}
\usepackage{amssymb}
\usepackage{subfig}
\usepackage{algorithm}
\usepackage{algpseudocode}

\newcommand{\R}{\mathbb{R}}
\newcommand{\ds}{\displaystyle}
\newcommand{\Id}{\textrm{Id}}
\newcommand{\x}{{\bf x}}
\newcommand{\Div}{{\rm div}}

\newtheorem{Theorem}{Theorem}[section]
\newtheorem{Lemma}{Lemma}[section]

\newtheorem{Proposition}{Proposition}[section]
\newtheorem{Corollary}{Corollary}[section]
\newtheorem{remark}{Remark}[section]

\newtheorem*{Assumption*}{Assumption}

\newtheorem{problem}{Problem}[section]
\newtheorem*{problem*}{Problem}
\numberwithin{equation}{section}

\begin{document}

\title{A new algorithm to determine the creation or depletion term of parabolic equations from boundary measurements}

\author{ Loc H. Nguyen 
\thanks{%
Department of Mathematics and Statistics, University of North Carolina at
Charlotte, Charlotte, NC 28223, USA, \text{loc.nguyen@uncc.edu}.}}

\date{}
\maketitle
\begin{abstract}
We propose a  robust  numerical method to find the coefficient of the creation or depletion term of parabolic equations from the measurement of the lateral Cauchy information of their solutions.
	Most papers in the field study this nonlinear and severely ill-posed problem using optimal control.
	The main drawback of this widely used approach is the need of some advanced knowledge of the true solution.
	In this paper, we propose a new method that opens a door to solve nonlinear inverse problems for parabolic equations without any initial guess of the true coefficient.
	This claim is confirmed numerically.
	The key point of the method is to derive a system of nonlinear elliptic equations for the Fourier coefficients of the solution to the governing equation with respect to  a special basis of $L^2$.
	We then solve this system by a predictor-corrector process, in which our computation to obtain the first and second predictors is effective.
	The  desired solution to the inverse problem under consideration follows.
\end{abstract}



\section{Introduction} \label{sec Intro}

Let $\Omega$ be a cube $(-R, R)^d \subset \R^d$, $d \geq 2$, where $R$ is a positive number.
Introduce a $d \times d$ matrix valued function $A$ with entries in the class $C^\infty(\R^d, \R^{d \times d})$. 
Assume throughout the paper that
\begin{enumerate}
    \item $A$ is symmetric; i.e, $A^T = A$,
    \item the matrix $A$ is uniformly elliptic; i.e., there exists a positive number $\mu$ such that 
    \[
        A(\x) \xi \cdot \xi \geq \mu |\xi|^2 \quad \mbox{for all } \x \in \R^d, \xi \in \R^d,
    \]
    \item for all $\x \in \R^d \setminus \Omega$, $A(\x) = \Id$ where $\Id$ is the identity matrix.
\end{enumerate}
Let ${\bf b}$ be a $d$-dimensional vector valued function in $C^\infty(\R^d, \R^d)$.
Define the operator
\begin{equation}
	Lw(\x) = \Div(A(\x) \nabla w(\x)) + {\bf b}{(\x)} \cdot \nabla w(\x) \quad \mbox{for all } w \in H^{2}(\R^d),\x \in \R^d.
	\label{L def}
\end{equation}
Consider the solution $u$ to the following initial value problem for the following parabolic equation
\begin{equation}
    \left\{
        \begin{array}{rcll}
             u_t(\x, t) &=& Lu(\x, t) + c(\x) u(\x, t) & \x \in \R^d, t \in (0, \infty),\\
             u(\x, 0) &=& f(\x) & \x \in \R^d.
        \end{array}
    \right.
    \label{para eqn}
\end{equation}
Here, $f$ is a smooth function defined on $\R^d$.
We refer the reader to \cite[Chapter 7]{Evans:PDEs2010} and \cite{LadyZhenskaya:ams1968} for the well-posedness of \eqref{para eqn}.
The regularity of the function $u(\x, t)$ can be found in those books. 
The second order term $\Div (A(\x)\nabla u(\x, t))$ describes the diffusion, the first order term ${\bf b}(\x) \cdot \nabla u(\x, t)$ describes the transport and the zero$^{\rm th}$ order term $c(\x) u(\x, t)$ describes creation or depletion. 
In this paper, we numerically solve the problem of reconstructing the coefficient $c(\x)$ of the creation or depletion term.
{
Roughly speaking, the creation or depletion term refers to the ability of produce or destroy, respectively, photons.
For e.g., in a chemical reaction, in which $u(\x, t)$ is the concentration of unstable gas,
the rate of decomposition is proportional to the concentration. This leads to the presence of $c(\x) u(\x, t)$ in \eqref{para eqn}. 
We would like to refer the reader to \cite[Chapter 6]{TikhonovSamarskii:1963} for details about this term.}
In the current paper, we propose a method to solve the following highly nonlinear and severely ill-posed coefficient inverse problem.

\begin{problem}	[Coefficient Inverse Problem (CIP)]
    Assume that $f(\x) \not = 0$ for all $\x \in \Omega$.
    Given a time $T > 0$ and the lateral Neumann data
    \begin{equation}
       F(\x, t) = u(\x, t) \quad \mbox{and } \quad  G(\x, t) = { A(\x) \nabla u(\x, t) \cdot \nu(\x)}
    \end{equation}
    for all $\x \in \partial \Omega$ and $t \in [0, T],$ determine the coefficient $c(\x)$, $\x \in \Omega$. 
    Here $\nu(\x)$ is the outward normal vector of $\partial \Omega$ at a point $\x$.
    \label{cip}
\end{problem}

Problem \ref{cip} has uncountable {practical} applications.
In fact, suppose that interior points of a medium are not accessible. 
In this case, by  measuring some boundary information of the function $u$, which are the heat and the heat flux in this paper, for a certain period of time  and by solving Problem \ref{cip}, one can determine the coefficient $c(\x)$ of the governing equation in \eqref{para eqn}, 
which enables us to inspect that medium without destructing it.
We recall here a specific example in bioheat transfer. 
In this field, the coefficient $c(\x)$ represents the blood perfusion. 
The knowledge of this coefficient plays a crucial role in calculating the temperature of the blood flowing through the tissue, see \cite{CaoLesnic:amm2019}.
{The uniqueness of Problem \ref{cip} is still open and is studied in an approximation context  of this paper.}
One can find the uniqueness of other versions of Problem \ref{cip} in \cite{BeilinaKlibanovBook, BukhgeimKlibanov:smd1981, PrilepkoKostin:RASBM1993} when some internal data are assumed to be known. When the Dirichlet to Neuman map is given, the reader can find the uniqueness in \cite{Isakov:ip1999}.
The uniqueness of Problem \ref{cip} is an assumption in this paper.
Another related problem is the inverse problem of recovery other coefficients, for e.g., the diffusion, or the initial condition from the measurement of the final time for parabolic equations. This problem is very important and interesting, see  \cite{klibanovYagola:arxiv2019, LiNguyen:IPSE2019, TuanKhoaAu:SIAM2019, Prilepko:pam2000,   Tuan:ip2017} for theoretical results and numerical methods. 
In this paper, we model the function $u$ defined on $\R^d \times [0, \infty)$.
However, the reader will see in Sections \ref{sec system} and \ref{sec method - quasi}, our analysis and algorithm are made inside $\Omega$. 
In other words, our method works when domain of the function $u$ and model \eqref{para eqn} is restricted to $\Omega$.

Coefficient inverse problems for parabolic equations were studied intensively.
Up to the knowledge of the author, the widely used method to solve this problem is the optimal control approach, see e.g., \cite{Borceaetal:ip2014, CaoLesnic:nmpde2018, CaoLesnic:amm2019, KeungZou:ip1998, YangYuDeng:amm2008} and references therein. 
The authors of \cite{Borceaetal:ip2014} applied the optimal control method involving a preconditioner to numerically compute the {thermal} conductivity with high quality.
The main drawback of this method is the need of a good initial guess for the true solution {which} is not always available. 
On the other hand, we specially draw the reader's attention to the convexification method, see \cite{KlibanovNik:ra2017, Klibanov:ip2015}, which can overcome the difficulty about the availability of the initial guess. 
In those papers \cite{KlibanovNik:ra2017, Klibanov:ip2015}, the authors introduce a convex functional whose minimizer yields the solution of the problem under consideration, by combining the quasi-reversibility method and the Carleman weight functions.
{One dimensional} numerical examples, illustrating the role of Carleman weight functions in convexifying the cost functionals,  are presented in \cite{KlibanovNik:ra2017}.
It is valuable to numerically test this convexification method in higher dimensions. 
We also cite to \cite{Nguyens:jiip2019} for another method to solve Problem \ref{cip} by repeatedly solving its linearization.
In the current paper, we propose a novel method in which no  advanced knowledge about the true coefficient is required.

%
Our method to solve Problem \ref{cip} consists of two stages. 
In the first stage, we eliminate the function $c(\x)$ from \eqref{para eqn}.
The resulting equation obtained in this stage is not a standard equation.
A numerical method to solve it is not available yet. 
In the second stage,
we approximate that nonstandard equation {as} a coupled system of elliptic partial differential equations. 
This system is derived  based on a truncation of the Fourier series, with respect to a special basis originally introduced in \cite{Klibanov:jiip2017}.
We apply a predictor-corrector procedure, in which the first approximation of the true solution is computed without any of its advance knowledge.
The solution of Problem \ref{cip} follows.

Two important steps in our method require us to find vector valued functions satisfying a system of elliptic partial differential equations and both Dirichlet and Neumann boundary conditions.
We employ the quasi-reversibility method, so-called a global minimization, for this purpose and we also prove the convergence of the quasi-reversibility method in our context, using a new Carleman estimate in \cite{NguyenLiKlibanov:IPI2019}.
The quasi-reversibility method was first introduced by Latt\`es and Lions in \cite{LattesLions:e1969} for numerical solutions of ill-posed problems for
partial differential equations. 
It has been studied intensively since then,
see e.g., \cite{Becacheelal:AIMS2015, Bourgeois:ip2006,
BourgeoisDarde:ip2010, BourgeoisPonomarevDarde:ipi2019, ClasonKlibanov:sjsc2007, Dadre:ipi2016, KaltenbacherRundell:ipi2019,
KlibanovSantosa:SIAMJAM1991, Klibanov:jiipp2013, LocNguyen:ip2019, NguyenLiKlibanov:IPI2019}. A survey
on this method can be found in \cite{Klibanov:anm2015}.

In Section \ref{sec system}, we derive the system mentioned above. 
In Section \ref{sec unique}, we study the uniqueness of this system. 
In Section \ref{sec method - quasi}, we propose  a numerical method to solve that system. Also in Section \ref{sec method - quasi}, we study the quasi-reversibility method that can be applied in our context.
In Section \ref{sec 4}, we describe the implementation using the finite difference method.
In Section \ref{sec numerical example}, we present some numerical results.
Section \ref{sec concluding} is for the concluding remarks.

\section{A nonlinear coupled system of elliptic equations} \label{sec system}

In this paper, the initial condition of $f$ is assumed to be {strictly positive in $\overline \Omega$} and in $H^{2 + \beta}(\R^d)$ for some $\beta > 0$ and has compact support in $\R^d$.
Then, since all coefficients of the operator $L$ are in class $C^{\infty}$, \eqref{para eqn} has a unique solution with $|u(\x, t)| \leq M$ and
 $u \in H^{2 + \beta, 1 + \beta/2} (\R^d \times [0, T])$ for some constant $M > 0$. 
 These unique solvability and regularity properties can be obtained by applying Theorem 6.1 in \cite[Chapter 5, \S 6]{LadyZhenskaya:ams1968} and Theorem 2.1 in \cite[Chapter 5, \S 2]{LadyZhenskaya:ams1968}.
We impose the condition that $\beta$ is large such that the function $u$ has second derivative with respect to time.
Further more, if we impose a stronger condition that $f$ belongs to the class $C_0^{\infty}(\R^d)$, then applying the induction arguments in the third paragraph of \cite[page 456]{LadyZhenskaya:ams1968}, we see that the function  is $k$ times differentiable  for any $k > 0$.

From now on, we denote by $\Omega_T$ the set $\Omega \times [0, T].$
Define the function 
\begin{equation}
	v(\x, t) = u_t(\x, t) \quad \mbox{for all } (\x, t) \in \Omega_T.
\end{equation}
It follows from \eqref{para eqn} that
\begin{equation}
	v_t(\x, t) = L v(\x, t)+ c(\x) v(\x, t) 
	\quad
	\mbox{for all } (\x, t) \in \Omega_T.
	\label{2.2}
\end{equation}
On the other hand, for all $\x \in \Omega,$
\begin{equation*}
	v(\x, 0) = u_t(\x, 0) = L f(\x) + c(\x) f(\x).
\end{equation*}
Therefore, 
\begin{equation}
	c(\x) = \frac{v(\x, 0) - L f(\x)}{f(\x)} \quad \mbox{for all } \x \in \Omega.
	\label{2.3}
\end{equation}
Plugging \eqref{2.3} into \eqref{2.2}, we obtain the following equation 
\begin{equation}
	v_t(\x, t) = Lv(\x, t) - \frac{Lf(\x)}{f(\x)} v(\x, t) + \frac{v(\x, 0) }{f(\x)} v(\x, t) 
	\quad {\mbox{for all } (\x, t) \in \Omega_T.}
	\label{v eqn}
\end{equation}
\begin{remark}
Solving the nonlinear equation \eqref{v eqn} is challenging due to the presence of the initial condition $v(\x, 0)$. 
A theoretical result to solve it is not yet available.
We employ the technique of truncating the Fourier series, see \cite{Klibanov:jiip2017}, to solve \eqref{v eqn}.
\end{remark}

Recall a special orthonormal  basis of $L^2(0, T)$ originally introduced by Klibanov \cite{Klibanov:jiip2017} in 2017. 
This basis plays a crucial role in deriving an approximate model whose solution will be used to directly compute the solution of Problem \ref{cip}.
For each $n \geq 1$, define the function $\phi_n(t) = (t - T/2)^{n - 1} \exp(t - T/2)$.
It is well-known that the set $\{\phi_n\}_{n = 1}^{\infty}$ is complete in $L^2(0, T)$. 
Employing the Gram-Schmidt orthonormalization process on this set, we obtain an orthonormal basis of $L^2(0, T)$. We denote this basis by $\{\Psi_n\}_{n = 1}^{\infty}.$
We have the {following} proposition.
\begin{Proposition}[see \cite{Klibanov:jiip2017}]
The basis $\{\Psi_n\}_{n = 1}^{\infty}$ satisfies the following properties:
\begin{enumerate}
	\item $\Psi_n$ is not identically zero for all $n \geq 1$,
	\item For all $m, n \geq 1$
	\[
		s_{mn} = \int_0^T \Psi_n'(t) \Psi_m(t) dt 
		= \left\{
			\begin{array}{ll}
				1 & \mbox{if } m = n,\\
				0 & \mbox{if } n < m.
			\end{array}
		\right.
	\] 
	As a result, for all integer $N > 0$, the matrix $S = (s_{mn})_{m, n = 1}^N$, is invertible.
\end{enumerate}
\label{prop MK}
\end{Proposition}

Recall the Fourier coefficients of the function $v(\x, t)$ 
\begin{equation}
	v_n(\x) = \int_{0}^T v(\x, t) \Psi_n(t) dt \quad \mbox{for all } \x \in \Omega.
	\label{2.5}
\end{equation}
We have
\[
	v(\x, t) = \sum_{n = 1}^{\infty} v_n(\x) \Psi_n(t) \quad \mbox{for all } (\x, t) \in \Omega_T.
\]
Fix a number $N > 0$. 
We approximate the function $v(\x, t)$ by the partial sum
\begin{equation}
	v(\x, t) = \sum_{n = 1}^{N} v_n(\x) \Psi_n(t) \quad \mbox{for all } (\x, t) \in \Omega_T.
	\label{2.6}
\end{equation}
In this approximation context, 
\begin{equation}
	v(\x, 0) = \sum_{n = 1}^{N} v_n(\x) \Psi_n(0) \quad \mbox{for all } \x \in \Omega
	\label{2.7}
\end{equation}
and
\begin{equation}
	v_t(\x, t) = \sum_{n = 1}^{N} v_n(\x) \Psi_n'(t) \quad \mbox{for all } {(\x, t)} \in \Omega_T.
	\label{2.8}
\end{equation}
Plugging \eqref{2.6}, \eqref{2.7} and \eqref{2.8} into \eqref{v eqn}, we have
\begin{equation*}
	\sum_{n = 1}^{N} v_n(\x) \Psi_n'(t) 
	= \sum_{n = 1}^{N} L v_n(\x) \Psi_n(t) 
	- \frac{Lf(\x)}{f(\x)} \sum_{n = 1}^{N} v_n(\x) \Psi_n(t) 
	+ \sum_{n, l = 1}^{N} \frac{ v_n(\x) \Psi_n(0) }{f(\x)}  v_l(\x) \Psi_l(t)
\end{equation*}
for all $(\x, t) \in \Omega_T.$
For each $m$ in $\{1, 2, \dots, N\}$, multiply both sides of the equation above by $\Psi_m(t)$ and then integrate the resulting equation with respect to $t$. 
Noting that
\[
	\int_0^T \Psi_n(t)\Psi_m(t) dt = \delta_{m - n},
\]
we have for each $m \in \{1, \dots, N\}$,
\begin{equation}
	L v_m(\x) - \sum_{n = 1}^N s_{mn} v_{n}(\x) - \frac{Lf(\x)}{f(\x)} v_m(\x) 
	+ \sum_{n = 1}^N \frac{\Psi_n(0)}{f(\x)}v_n(\x) v_m(\x) = 0 
	\label{2.9}
\end{equation}
for all $\x \in \Omega.$
On the other hand, for each $m \in \{1, \dots, N\}$, the function $v_m(\x)$ satisfies the following constraints for all $\x \in \partial \Omega$
{
\begin{equation}
	\begin{array}{c}
	\ds v_m(\x) = F_m(\x) = \int_0^T F_t(\x, t) \Psi_m(t) dt, 
\\
	 \ds \partial_\nu v_m(\x) = G_m(\x) = \int_0^T G_t(\x, t) \Psi_m(t) dt,
	\end{array}
	\label{2.10}
\end{equation}

\begin{remark}[Data]
	From now on, we consider $F_m$ and $G_m$ as the {Cauchy} indirect data. These data are computed directly from the time dependent data via \eqref{2.10}.
	\end{remark}
	\begin{remark}[Noise]
	Let $F^*$ and $G^*$ be the data without noise. Let ${\rm rand}$ be the function of uniformly distributed random numbers in the range $[-1, 1]$. 
	For $\delta > 0$, define
	\[
		F^\delta = F^* + \delta{\rm rand} \quad G^\delta = G^* + \delta {\rm rand}.
	\]
	The corresponding noisy indirect data are denoted by $F^\delta_m$ and $G^\delta_m$ respectively.
	There is a difficulty in computing $F^\delta_m$ and $G^\delta_m$ via \eqref{2.10}.
	Computing the derivatives $F_t$ and $G_t$ of the data  for \eqref{2.10}  is unstable.
	 We use the well-known Tikhonov regularization technique to compute them.	
	In this paper, we test our numerical method for simulated data with noise level $\delta = 5\%$.
	\label{rem noise}
\end{remark}

In summary, we have proved the following proposition.
\begin{Proposition}
	Fix $N > 0$.
	Assume that the function $v(\x, t)$, {for} $(\x, t) \in \Omega_T$, can be well-approximated by the expression in \eqref{2.6} with the function $v_n(\x)$, $n \in \{1, 2, \dots, N\}$, given in \eqref{2.5}. 
	Then the Fourier coefficients $v_n$ satisfies the over-determined system of partial differential equations \eqref{2.9}--\eqref{2.10}.
	\label{prop 2.2}
\end{Proposition}

\section{The uniqueness of Problem \ref{cip} in the approximation context} \label{sec unique}

Proposition \ref{prop 2.2} suggests a method to numerically solve Problem \ref{cip}.
We solve the nonlinear system \eqref{2.9}--\eqref{2.10} for a vector $(v_1, \dots, v_N)$ and compute the coefficient $c$ via \eqref{2.6} and  then \eqref{2.3}. 
Hence, the uniqueness of the solution to the nonlinear system \eqref{2.9}--\eqref{2.10} implies the unique reconstruction of Problem \ref{cip} assuming the approximation \eqref{2.6}. 
In this section, we establish the uniqueness for \eqref{2.9}--\eqref{2.10}.
It is sufficient to study this uniqueness for the bounded solution.
In fact, due to the conditions imposed on the source function $f$ in the first paragraph of Section \ref{sec system}, namely $f \in C_0^{\infty}(\R^d)$, the true function $u(\x, t)$ is $k$ times differentiable for any $k > 0$. 
Thus, the function $v(\x, t) = u_t(\x, t)$ is bounded. 
It follows from \eqref{2.5} that the true solution to \eqref{2.9}--\eqref{2.10} is in $L^{\infty}(\Omega)$.

In order to establish the uniqueness of the bounded solution to \eqref{2.9}--\eqref{2.10}, we need the following Carleman estimate.
\begin{Lemma}[Carleman estimate]
Let the number $b>R$. 
Then there exist numbers $p_0 \geq 1$ and $\lambda_0 \geq 1$ depending only on $\mu$, $b$, $d$, $R$, $\|A\|_{L^{\infty}(\Omega)^{d \times d}}$
 such that the following Carleman estimate holds:
\begin{equation}
\int_{\Omega } |\Div(A \nabla u)|^{2}\exp \left[ 2\lambda \left(
x_d+b\right) ^{p}\right] d\x
\geq 
C \lambda \int_{\Omega }
\left[ |\nabla u|^{2}+\lambda^{2} u^{2}\right] \exp \left[
2\lambda \left( x_d+b\right) ^{p}\right] d\mathbf{x},  
\label{3.222222}
\end{equation}%
for all $\lambda \geq \lambda _{0},$ $p\geq p_{0}$ and $u \in
H^2(\Omega)$ with $u = \partial_{\nu} 0 = 0$ on $\partial \Omega$. 
Here, the constant $C$ depends only on $\mu$, $b$, $d$, $R$ and $\|A\|_{L^{\infty}(\Omega)^{d \times d}}$. 
\label{lem Carleman}
\end{Lemma}

Lemma \ref{lem Carleman} is a direct consequence of \cite[Theorem 4.1]{NguyenLiKlibanov:IPI2019}.
We do not repeat the proof in this paper. 
This Lemma plays a crucial role in the proof of Theorem \ref{thm uniqueness} below. Moreover, it will be applied to prove the convergence of the quasi-reversibility method, see Theorem \ref{thm}.

We have the theorem.

\begin{Theorem}
	The nonlinear system \eqref{2.9}--\eqref{2.10} has at most one solution in $H^2(\Omega) \cap L^{\infty}(\Omega)$.
	\label{thm uniqueness}
\end{Theorem}
\begin{proof}	 
	Let $U = (u_1, \dots, u_N)$ and $V = (v_1, \dots, v_N)$ be two solutions in $H^2(\Omega) \cap L^{\infty}(\Omega)$ to \eqref{2.9}--\eqref{2.10}.
	It follows from \eqref{2.9} that for all $\x \in \Omega$ and $m \in \{1, \dots, N\}$,
\begin{equation}
	L u_m(\x) - \sum_{n = 1}^N s_{mn} u_{n}(\x) - \frac{Lf(\x)}{f(\x)} u_m(\x) 
	+ \sum_{n = 1}^N \frac{\Psi_n(0)}{f(\x)}u_n(\x) u_m(\x) = 0 
	\label{13}
\end{equation}
and
\begin{equation}
	L v_m(\x) - \sum_{n = 1}^N s_{mn} v_{n}(\x) - \frac{Lf(\x)}{f(\x)} v_m(\x) 
	+ \sum_{n = 1}^N \frac{\Psi_n(0)}{f(\x)}v_n(\x) v_m(\x) = 0 
	\label{14}
\end{equation}
It follows from \eqref{13} and \eqref{14} that for all $\x \in \Omega$ and $m \in \{1, \dots, N\}$,
\begin{multline}
	L (u_m(\x) - v_m(\x)) - \sum_{n = 1}^N s_{mn} (u_n(\x) - v_n(\x)) - \frac{Lf(\x)}{f(\x)} (u_m(\x) - v_m(\x))
	\\
	+ \sum_{n = 1}^N \frac{\Psi_n(0)}{f(\x)}[u_n(\x) (u_m(\x) - v_m(\x)) + v_m(\x) (u_n(\x) - v_n(\x))]
	= 0.
	\label{15}
\end{multline}
Let $h = (h_1, \dots, h_N) = U - V$.
Since $U$ and $V$ are in $L^{\infty}(\Omega),$ we can find a number $M$ such that $|u_m| \leq M$ and $|v_m| \leq M$ for all $m \in \{1, \dots, N\}.$  
It follows from \eqref{L def}, \eqref{2.10} and \eqref{15} that 
\begin{equation}
	\left\{
		\begin{array}{rcll}
			\ds\sum_{m = 1}^N|\Div(A \nabla h_m)|^2 &\leq& C_1\Big(\ds\sum_{m = 1}^N |h_m|^2 + \sum_{m = 1}^N|\nabla h_m|^2\Big) &\mbox{in } \Omega,\\
			h_m &=& 0 &\mbox{on } \partial \Omega,\\
			\partial_{\nu} h_m &=& 0 &\mbox{on } \partial \Omega
		\end{array}
	\right. 
	\label{18}
\end{equation}
for some constant $C_1$ depending only on $N$, $M$, $(s_{mn})_{m, n = 1}^N$, $(\Psi_m)_{m = 1}^N$, $f$, and ${\bf b}$.
Let $b > 0$, $\lambda > \lambda_0$ and $p > p_0$ where $\lambda_0$ and $p_0$ are in Lemma \ref{lem Carleman}.
Applying the Carleman estimate in Lemma \ref{lem Carleman} for each function $h_m$, $m \in {1, \dots, N}$, we have
\begin{equation}
\sum_{m = 1}^N\int_{\Omega } |\Div(A \nabla h_m)|^{2}\exp [ 2\lambda (
x_d+b)^{p}] d\x
\geq 
C \lambda \sum_{m = 1}^N  \int_{\Omega }
\left[ |\nabla h_m|^{2}+\lambda^{2} |h_m|^{2}\right] \exp \left[
2\lambda \left( x_d+b\right) ^{p}\right] d\mathbf{x}.
\label{19}
\end{equation}
It follows from \eqref{18} and \eqref{19} that
\begin{multline*}
	C_1\int_{\Omega}\Big(\ds\sum_{m = 1}^N |h_m|^2 + \sum_{m = 1}^N|\nabla h_m|^2\Big) \exp \left[
2\lambda \left( x_d+b\right) ^{p}\right]d\x 
\\
\geq 
	C \lambda \sum_{m = 1}^N  \int_{\Omega }
\left[ |\nabla h_m|^{2}+\lambda^{2} |h_m|^{2}\right] \exp \left[
2\lambda \left( x_d+b\right) ^{p}\right] d\mathbf{x}.
\end{multline*}
Choosing $\lambda$ sufficiently large, we obtain $h_m = 0$ for all $m \in \{1, \dots, N\}.$
\end{proof}

\begin{remark}	
	The uniqueness of the solution $(v_1, \dots, v_N)$ to the nonlinear system \eqref{2.9}--\eqref{2.10} implies the unique reconstruction of the solution to Problem \ref{cip}. 
	This can be seen via the reconstruction formulas \eqref{2.6} and \eqref{2.3}.
	This uniqueness holds true only in the approximation context \eqref{2.6} while the uniqueness for the true model in the time domain is extremely challenging, which is out of the scope of this paper.
\end{remark}

\section{The method to solve Problem \ref{cip}} \label{sec method - quasi}

Solving Problem \ref{cip} is reduced to {solving} the system of nonlinear partial differential {equations} \eqref{2.9}--\eqref{2.10}.

\subsection{An iterative process} \label{sec method}

We propose the following iterative method, in which a predictor-corrector procedure is applied.
The first predictor, named as $V^{(0)}$, is set to be the solution of the linear system obtained by removing from \eqref{2.9} the nonlinear term. 
More precisely, we set $V^{(0)} = (v_1^{(0)}, \dots, v_N^{0})^T$ as the solution of 
\begin{equation}
	L v_m^{(0)}(\x) - \ds\sum_{n = 1}^N s_{mn} v_{n}^{(0)}(\x) - \frac{Lf(\x)}{f(\x)} v_m^{(0)}(\x)  = 0 
	\quad \mbox{for all } \x \in \Omega
	\label{3.1}
\end{equation}
and
\begin{equation}
	v_m^{(0)}(\x) = F_m(\x), 
	\quad
	\partial_{\nu}v_m^{(0)}(\x) = G_m(\x) 
	\quad \mbox{for all } \x \in \partial \Omega
	\label{3.2}
\end{equation}
for $m \in \{1, \dots, N\}.$
Next, by induction, assume that $V^{(p)}$ is known for some positive integer $p$, we find $V^{(p + 1)}$ by solving the equation obtained from \eqref{2.9} by replacing $v_m$ in the
nonlinear term by its approximation $v_m^{(p)}.$ 
That means, $V^{(p + 1)} = (v_1^{(p + 1)}, \dots, v_N^{(p + 1)})^T$ is set to be the solution of 
\begin{equation}
	L v_m^{(p + 1)}(\x) - \sum_{n = 1}^N s_{mn} v_{n}^{(p + 1)}(\x) - \frac{Lf(\x)}{f(\x)} v_m^{(p + 1)}(\x) 
	+ \sum_{n = 1}^N \frac{\Psi_n(0)}{f(\x)}v_n^{(p+1)}(\x) v_m^{(p)}(\x) = 0 \quad \mbox{for all } \x \in \Omega
	\label{3.3}
\end{equation}
 and
\begin{equation}
	v_m^{(p + 1)}(\x) = F_m(\x), \
	\partial_{\nu} 	v_m^{(p + 1)}(\x) = G_m(\x) \quad \mbox{for all } \x \in \partial \Omega
	\label{3.4}
\end{equation}
for each $m \in \{1, \dots, N\}.$

Due to the presence of the latteral Cauchy data, both problems \eqref{3.1}--\eqref{3.2} and \eqref{3.3}--\eqref{3.4} are over-determined.
In the case when the data is noisy, they might not have a solution. 
We, therefore, employ the quasi-reversibility method to solve them. In Proposition \ref{Prop 3.1}, we show the existence of an ``approximation" of the true solution in the case when the data has no noise. In Theorem \ref{thm}, we prove the convergence of this approximation as the noise tends to $0$.

\subsection{The quasi-reversibility method}\label{sec 3.2}

We next recall the quasi-reversibility method to solve systems of partial differential equations with Cauchy boundary data.
The two systems of elliptic partial differential equations \eqref{3.1}--\eqref{3.2} and \eqref{3.3}--\eqref{3.4} are over-determined due to both Dirichlet and Neumann boundary conditions imposed.
We use the quasi-reversibility method to solve them.
A general form of the {problem \eqref{3.1}--\eqref{3.2} and the problem \eqref{3.3}--\eqref{3.4}} is given by
\begin{equation}
	\left\{
	\begin{array}{rcll}
		\Div(A \nabla V) + BV &=& 0 &\x \in \Omega,\\
		V &=& \mathcal{F} & \x \in \partial \Omega,\\
		\partial_{\nu}V &=& \mathcal{G} & \x \in \partial \Omega
	\end{array}
	\right.
	\label{3.5}
\end{equation}
where $A$ is introduced in Section \ref{sec Intro} and $B$ is a $N \times N$ matrix valued function in $L^{\infty}(\Omega).$
We have the {following} proposition whose proof closely follows that of Theorem 3.1 in \cite{NguyenLiKlibanov:IPI2019}. 
\begin{Proposition} 	
	Fix $\epsilon > 0$. 
	Then, the functional 
	\begin{equation}
		J_{\epsilon}(V) = \int_{\Omega}|\Div(A \nabla V) + BV|^2 d\x 
		+ \int_{\partial \Omega} |V - \mathcal{F}|^2 d\sigma(\x) 
		+ \int_{\partial \Omega} |\partial_{\nu} V - \mathcal{G}|^2 d\sigma(\x)
		+ \epsilon \|V\|^2_{H^2(\Omega)}
		\label{3.7}
	\end{equation}
	has a unique minimizer on $H^2(\Omega).$
	This minimizer $V_{\epsilon}$ is called the regularized solution of \eqref{3.5}.
	\label{Prop 3.1}
\end{Proposition}

In summary, we propose Algorithm \ref{alg} to solve Problem \ref{cip} via solving \eqref{2.9}--\eqref{2.10} by the quasi-reversibility method.

\begin{Proposition}[The uniqueness of \eqref{3.5}]
The system \eqref{3.5} has at most one solution in $H^2(\Omega)$.
\label{Prop unique}
\end{Proposition}
The proof of this proposition is very similar to that of Theorem \ref{thm uniqueness} by using the Carleman estimate. We do not repeat the proof here.

\begin{Theorem}[The convergence of the quasi-reversibility method for \eqref{3.5}]
{Assume that there uniquely exists a true solution $V^*$ to \eqref{3.5} with the boundary data $\mathcal{F}$ and $\mathcal{G}$ replaced by the corresponding noiseless ones, denoted by $\mathcal{F}^*$ and $\mathcal{G}^*$ respectively.} 
Let $\mathcal{F}^{\delta}$ and $\mathcal{G}^{\delta}$ be the corresponding noisy data for some $\delta > 0$.
Assume that there exists an ``error" vector valued function $\mathcal{E}$ such that
\begin{equation}
	\mathcal{E} = \mathcal{F}^{\delta} - \mathcal{F}^{*} \quad \mbox{and } \quad \partial_{\nu} \mathcal{E} = \mathcal{G}^{\delta} - \mathcal{G}^{*} \quad {\mbox{on } \partial \Omega}
	\label{error boundary}
\end{equation} 
and assume that
\begin{equation}
	\|\mathcal{E}\|_{H^2(\Omega)} \leq \delta.
	\label{error est}
\end{equation}
Then, $V_{\epsilon}^{\delta}$, {the minimizer of $J_{\epsilon}$ with $\mathcal F$ and $\mathcal G$ replaced by $\mathcal F^{\delta}$ and $\mathcal G^{\delta}$ respectively}, satisfies the estimate
\begin{equation}
 \|V^{\delta}_{\epsilon} - V^*\|^2_{H^1(\Omega)^N} \leq C(\delta^2 + \epsilon \|V^*\|_{H^2(\Omega)^N}^2).
 \label{23}
\end{equation}
\label{thm}
\end{Theorem}

\begin{remark}
	The convergence for the quasi-reversibility method guaranteed by Theorem \ref{thm} is similar to that in \cite[
Theorem 5.1]{NguyenLiKlibanov:IPI2019}. 
The main difference of two results is in the objective functional is minimized subject to some boundary constraints while in the current paper, such constraints are relaxed by adding the two boundary integrals in \eqref{3.7}.
\end{remark}

\begin{remark}[The existence of the error function $\mathcal{E}$ in the statement of Theorem \ref{thm}]
	By employing the main result in \cite{BourgeoisDarde:ip2010},
	we can prove the existence of the vector valued function $\mathcal{E}$ that satisfies \eqref{error boundary} if we
	 impose the reasonable conditions that $\Omega$ is in the class $C^{1,1}$ and that $(\mathcal{F}^{\delta} - \mathcal{F}^{*}, \mathcal{G}^{\delta} - \mathcal{G}^{*})$ belongs to $H^{3/2}(\partial \Omega)$.
\end{remark}

\begin{proof}[Proof of Theorem \ref{thm}]
	Since $V_{\epsilon}^{\delta}$ is the regularized solution to \eqref{3.5}, it is the minimizer of $J_{\epsilon}$, defined in \eqref{3.7}.
	Hence, for all $\phi \in H^2(\Omega)^N$, we have
	\begin{multline}
		\langle \Div(A \nabla V_{\epsilon}^{\delta}) + BV_{\epsilon}^{\delta}, \Div(A \nabla \phi) + B \phi \rangle_{L^2(\Omega)^N}
		\\
		+ \langle V_{\epsilon}^{\delta} - \mathcal{F}^{\delta}, \phi\rangle_{L^2(\partial\Omega)^N}
		+\langle \partial_{\nu}V_{\epsilon}^{\delta} - \mathcal{G}^{\delta}, \partial_{\nu}\phi\rangle_{L^2(\partial\Omega)^N} 
		+ \epsilon \langle V_{\epsilon}^{\delta}, \phi\rangle_{H^2(\Omega)^N} = 0
		\label{3.9}
	\end{multline}
	for all $\phi \in H^2(\Omega)^N.$
	On the other hand, since $V^{*}$ is the true solution to \eqref{3.5}
	\begin{multline}
		\langle \Div(A \nabla V^*) + BV^*, \Div(A \nabla \phi) + B \phi \rangle_{L^2(\Omega)^N}
		+ \langle V^* - \mathcal{F}^*, \phi\rangle_{L^2(\partial\Omega)^N}
		\\
		+\langle \partial_{\nu}V^* - \mathcal{G}^*, \partial_{\nu}\phi\rangle_{L^2(\partial\Omega)^N} 
		+\epsilon \langle V^*, \phi\rangle_{H^2(\Omega)^N} 
		=
		 \epsilon \langle V^*, \phi\rangle_{H^2(\Omega)^N}
		 \label{3.10}
	\end{multline}
	for all $\phi \in H^2(\Omega)^N.$
	Taking the difference of \eqref{3.9} and \eqref{3.10}, we have
	\begin{multline}
		\langle \Div(A \nabla W) + BW,  \Div(A \nabla \phi) + B\phi\rangle_{L^2(\Omega)^N} 
		+\langle W - (\mathcal{F}^\delta - \mathcal{F}^*), \phi\rangle_{L^2(\partial\Omega)^N}\\
		+\langle \partial_{\nu}W - (\mathcal{G}^\delta - \mathcal{G}^*), \partial_{\nu}\phi\rangle_{L^2(\partial\Omega)^N}
		+ \epsilon \langle W, \phi\rangle_{H^1(\Omega)^N}
		= -\epsilon \langle V^*, \phi\rangle_{H^2(\Omega)^N}
		\label{3.1313}
	\end{multline}
where $W = V^{\delta}_{\epsilon} - V^*$ for all $\phi \in H^2(\Omega)^N.$
Using
\begin{equation}
	\phi = W - \mathcal{E} = V^{\delta}_{\epsilon} - V^* - \mathcal{E}
	\label{3.13}
\end{equation}
as a test function in \eqref{3.1313}
 and using \eqref{error boundary}, we have
\begin{multline*}
	\|\Div(A \nabla \phi) + B\phi\|_{L^2(\Omega)^N}^2
	+ \langle \Div(A \nabla \mathcal{E}) + B\mathcal{E},  \Div(A \nabla \phi) + B\phi\rangle_{L^2(\Omega)^N}
	\\
	+ \|\phi\|_{L^2(\partial\Omega)^N}^2
	+ \|\partial_{\nu}\phi\|_{L^2(\partial\Omega)^N}^2
	+ \epsilon \| \phi\|_{H^2(\Omega)^N}^2
	+ \epsilon \langle \mathcal{E}, \phi\rangle_{H^2(\Omega)^N} 
	= -\epsilon \langle V^*, \phi\rangle_{H^2(\Omega)^N}.
\end{multline*}
Applying the inequality $|\langle u, v\rangle| \leq 1/2(\|u\|^2 + \|v\|^2)$, \eqref{error est} and the trace theory, we have
\begin{equation}
	\|\Div(A \nabla \phi) + B\phi\|_{L^2(\Omega)^N}^2 \leq C(\delta^2 + \epsilon \|V^*\|_{H^2(\Omega)^N}^2).
	\label{step 1}
\end{equation}
Here, $C$ is a generic constant that might change from estimate to estimate.
Choose $b > R$, $\lambda > \lambda_0$, $p \geq p_0$ where $b,$ $\lambda_0$ and $p_0$ are as in Lemma \ref{lem Carleman}.
It is not hard to verify that the function $\phi$ satisfies the homogenous boundary conditions $\phi = \partial_{\nu} \phi = 0$ on $\partial_{\Omega} \times [0, T]$.
Using the Carleman estimate in Lemma \ref{lem Carleman}, we can bound the left hand side of \eqref{step 1} as follows
\begin{align*}
	\int_{\Omega} |\Div(&A \nabla \phi) + B\phi|^2 d\x
	\\
	&\geq 
	 \exp(-2\lambda (R + b))\int_{\Omega} \exp(2\lambda (x_d + b)) |\Div(A \nabla \phi) + B\phi|^2 d\x
	 \\
	 &\geq C\int_{\Omega} [\exp(2\lambda (x_d + b)) |\Div(A \nabla \phi)|^2 -  \exp(2\lambda (x_d + b)) |B\phi|^2] d\x
	 \\
	 &\geq C \int_{\Omega} [\exp(2\lambda (x_d + b)) (\lambda^3 |\phi|^2 + \lambda |\nabla \phi|^2) -  \exp(2\lambda (x_d + b)) |B\phi|^2] d\x.
\end{align*}
Choosing $\lambda$ sufficiently large, since $B \in L^{\infty}(\Omega),$ we have 
\[
	\|\Div(A \nabla \phi) + B\phi\|_{L^2(\Omega)^N}^2 \geq 
	C \|\phi\|^2_{H^1(\Omega)^N}.
\]
This, together with \eqref{3.13} and \eqref{step 1}, implies {\eqref{23}.
The theorem is proved.}
\end{proof}

\begin{Corollary}
	Theorem \ref{thm} implies that the regularized solution to \eqref{3.5} obtained by quasi-reversibility method is close to the true solution.  
	In fact, if $\epsilon < O(\delta^2)$, it follows from \eqref{23} that
\[
	\|V^{\delta}_{\epsilon} - V^*\|_{H^1(\Omega)^N} \leq C\delta \|V^*\|_{H^2(\Omega)^N}.
\]	The rate of convergence is Lipschitz.
\end{Corollary}

	Theorem \ref{thm} guarantees that Steps \ref{reg 1} and \ref{reg 2} in Algorithm \ref{alg} below provide good approximations of the sequence $\{c_p^{\delta}\}_{p = 1}^{\infty}$ in comparison to the sequence $\{c_p^*\}_{p = 1}^{\infty}$ with the Lipschitz rate provided that $\epsilon = O(\delta^2)$ as $\delta \to 0^+$. 
	If the sequence $\{c_p^*\}_{p = 1}^{\infty}$ converges to the solution of Problem \ref{cip}, Algorithm \ref{alg} yields a numerical procedure to solve it.
	This convergence is verified numerically in Section \ref{sec numerical example}.

\subsection{The procedure to solve the coefficient inverse problem for parabolic equations}

By Proposition \ref{prop 2.2}, the strategy to solve \eqref{2.9}--\eqref{2.10} described in Section \ref{sec method} and the convergence of the quasi-reversibility method, see Theorem \ref{thm},  we propose Algorithm \ref{alg} to reconstruct the coefficient $c(\x)$ {for} $\x \in \Omega.$
\begin{algorithm}
\caption{\label{alg}The procedure to solve Problem \ref{cip}}\label{euclid}
\begin{algorithmic}[1]
	\State\, \label{st1} Choose a number $N$. Construct the functions $\Psi_m$, $1 \leq m \leq N,$  and compute the matrix S as in Proposition \ref{prop MK}. Fix $\epsilon > 0.$
	\State\, \label{reg 1} Find the regularized solution  $V^{(0)}$  of \eqref{3.1}--\eqref{3.2}.
	\State\, \label{st3} Compute $v^{(0)}$ via \eqref{2.6} and $c^{(0)}$ via \eqref{2.3}.
	\State\, \label{reg 2} Assume that we know $V^{(p)}$ and $c^{(p)}$.  Set $V^{(p + 1)}$ as the regularized solution to \eqref{3.3}--\eqref{3.4}.
	\State\, \label{st5} Compute $v^{(p + 1)}$ via \eqref{2.6} and $c^{(p+1)}$ via \eqref{2.3}.
	\State\, \label{st6} Define
	\[
		{E}(p) = \frac{\|c^{(p)} - c^{(p + 1)}\|_{L^{\infty}(\Omega)}}{\|c^{(p + 1)}\|_{L^{\infty}(\Omega)}}.
	\] Choose $c = c^{(p^*)}$ for $p^*$ such that ${E}(p^*)$ is sufficiently small.
	\end{algorithmic}	
\end{algorithm}

\begin{remark}
	Unlike the widely used least squares method to solve ill-posed inverse problems, we do not require a good initial guess for the true coefficient $c(\x)$. 
	Our first approximation is computed in Steps \ref{reg 1} and \ref{st3} of Algorithm \ref{alg}. 
	It is shown in Section \ref{sec numerical example} that the functions $c^{(0)}$ are acceptable in Tests 1, 3 and 4.
	In contrast, our computed $c^{(0)}$ is poor in Test 2. However, the error is automatically corrected when we find $c^{(1)}$ in Steps \ref{reg 2} and  \ref{st5}.
\end{remark}

\section{The implementation using the finite difference method} \label{sec 4}

We test our method in the simple case when $d = 2$, $\Omega = (-1, 1)^2$, $L = \Delta$ is the Laplacian and $f(\x) = 100$.

\subsection{The forward problem}
To generate the simulated data, we solve the forward problem of Problem \ref{cip}. 
That means, given $c(\x)$ {(see each test below for the definition of $c(\x)$)}, we need to compute the solution $u(\x, t)$ to \eqref{para eqn} on the whole plane $\R^2$. 
Instead of doing so, we solve an analog of \eqref{para eqn} on a domain $\Omega_1 = ({-}R_1, R_1)^2$ where $R_1 = 3 > R = 1.$
{
This domain approximation does not effect our analysis because the formulation of  the parabolic equation in \eqref{para eqn} is used only inside the domain $\Omega$ when we study the inverse problem.
}

In other words, we solve the equation
\begin{equation}
	\left\{
		\begin{array}{rcll}
			u_t(\x, t) &=& \Delta u(\x, t) + c(\x) u(\x, t) &\x \in \Omega_1, \quad t \in [0, T],\\
			u(\x, t) &=& f(\x) &\x \in \partial \Omega_1, \quad t \in [0, T],\\
			u(\x, 0) &=& f(\x) &\x \in \Omega_1.
		\end{array}
	\right.
\label{eqn Omega1}
\end{equation}
Here, we choose the time-independent Dirichlet boundary data for the simplicity. 
In this paper, we solve problem \eqref{eqn Omega1} by implicit method using finite differences by {the backward Euler scheme}. 
In the finite difference scheme, we find the function $u(\x, t)$ on the grid of points 
\begin{multline}
	\Big\{(x_i = -R_1 + (i-1)d_{\x_1}, y_j = -R_1 + (j-1)d_{\x_1}, (l - 1)d_t): 
	\\
	1 \leq i, j \leq N_1, 1 \leq l \leq N_t\Big\} \subset \overline \Omega_1 \times [0, T],
\end{multline} where $N_1$ and $N_t$ are two large integers, $d_{\x_1} = 2R_1/(N_1 - 1)$ and $d_t = T/(N_t - 1)$.
In our computational program, $N_1$ is set to be 240 and $N_t = 100.$
Having the function $u(\x, t)$ for all $\x \in \Omega_1$ and $t \in [0, T]$ in hand, we can directly extract the data $F(\x, t) = u(\x, t)$ and $G(\x, t) = \partial_{\nu} u(\x, t)$ on $\partial \Omega \times [0, T]$.

\subsection{The inverse problem}
In this section, we present how to implement Algorithm \ref{alg} in the finite difference scheme.
Similarly to the previous section, we define and then compute the function $c$ on a uniform grid of points 
\begin{multline*}
	\{(x_i = -R + (i-1)d_\x, y_j = -R + (j-1)d_\x, (l - 1)d_t):
	\\ 1 \leq i, j \leq N_\x, 1 \leq l \leq N_t\} \subset \overline \Omega \times [0, T],
\end{multline*} where $N_\x$ and $N_t$ are two large integers, $d_\x = 2R/(N_\x - 1)$ and $d_t = T/(N_t - 1)$.
In all numerical tests in Section \ref{sec numerical example}, we take $N_\x = 80$ and $N_t = 100.$

We next present each step of Algorithm \ref{alg}.

\noindent {\bf Step \ref{st1}.} In this step, to choose ``truncation" number $N$.
To do so, we take a ``reference" function $v$ in one of the examples in Section \ref{sec numerical example} and then compute the absolute difference 
\[
	e_N(\x, t) = \Big|v(\x, t)  - \sum_{n = 1}^N v_n(\x)\Psi_n(t)\Big|
	\quad 
	\mbox{for all } \x \in \Omega, \quad t \in [0, T].
\] 
We observe that the larger $N$, the smaller $\|e_N\|_{L^{\infty}(\Omega \times [0, T])}$. 
We examine the function $e_N$ when $N = 5,$ $N = 10$ and $N = 25$, see Figure \ref{fig trunc}. 
It is evident from Figure \ref{fig 1 25} that when $N = 25$, $\|e_N\|_{L^{\infty}}$ is sufficiently small, about $4 {\times} 10^{-3}$. 
\begin{figure}
\begin{center}
	\subfloat[The function $e_5(\x, t = 0.3)$]{\includegraphics[width = 0.3\textwidth]{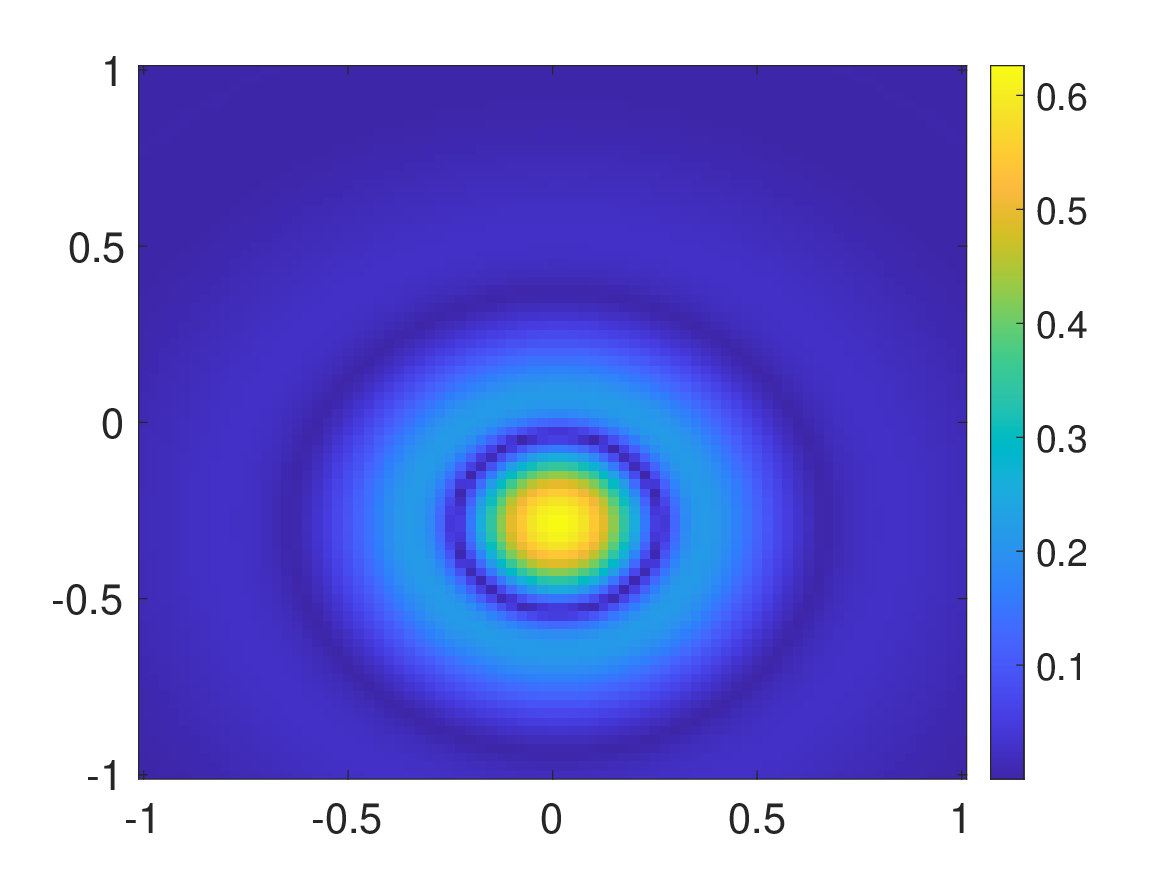}}
	\hfill
	\subfloat[The function $e_{10}(\x, t = 0.3)$]{\includegraphics[width = 0.3\textwidth]{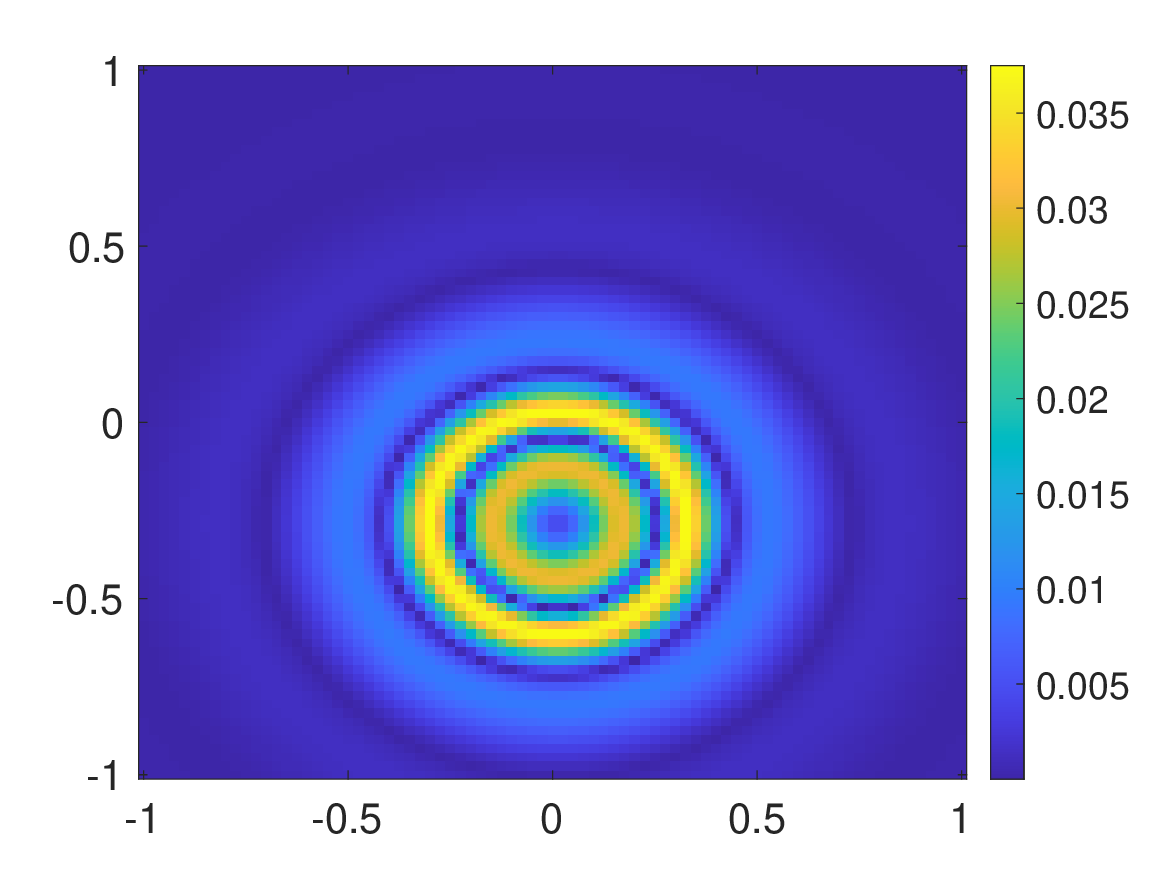}}
	\hfill
	\subfloat[\label{fig 1 25}The function $e_{25}(\x, t = 0.3)$]{\includegraphics[width = 0.3\textwidth]{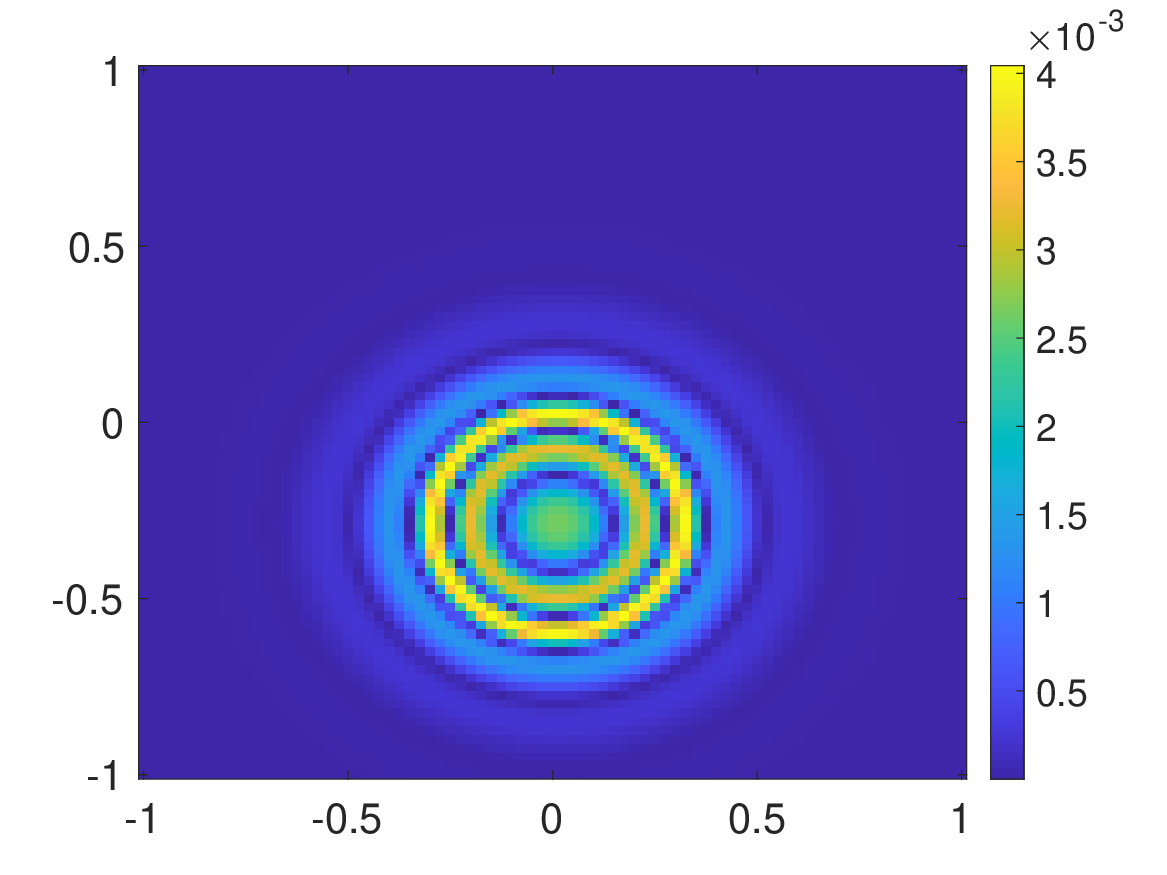}}
	\caption{\label{fig trunc}The difference of the function $v(\x, t)$ and the function $\sum_{n = 1}^N v_n(\x) \Psi_n(t)$ when $t = 0.3$.}
\end{center}
\end{figure}

As a result, we choose $N = 25$. We use this choice of $N$ for all numerical tests.
We observe that using higher $N$ does not improve the quality of the reconstructed coefficient $c(\x)$. 
Also in Step \ref{st1} of Algorithm \ref{alg}, we choose the regularized parameter $\epsilon = 10^{-9}.$ {This choice of $\epsilon$ is based on a trial-error process in the case when the given data is noisy.}

\noindent {\bf Step \ref{reg 1}.} Compute the vector valued function $V^{(0)}$, which is set to be the minimizer of the functional, due to the quasi-reversibility method,
\begin{multline}
	J_{\epsilon}^{(0)}(V) = \sum_{m = 1}^N\Big[\int_{\Omega} 
	\Big|\Delta v_m(\x) - \sum_{n = 1}^N s_{mn} v_{n}(\x) - \frac{\Delta f(\x)}{f(\x)} v_m(\x)\Big|^2 d\x
	\\
	+ \int_{\partial \Omega} |v_m(\x) - F_m(\x)|^2d\sigma(\x)
	+ \int_{\partial \Omega} |\partial_{\nu}v_m(\x) - G_m(\x)|^2d\sigma(\x)
	\\
	+ \epsilon \int_{\Omega} (|v_m(\x)|^2 + |\nabla v_m(\x)|^2) d\x
	\Big].
	\label{4.2}
\end{multline}
Here for each $m \in \{1, \dots, N\},$ $F_m = F_m^\delta$ and $G_m = G_m^\delta$ are computed in \eqref{2.10} with $F^\delta$ and $G^{\delta}$ replacing $F$ and $G$ respectively with $\delta = 5\%$.

Here, we replace the term $\Vert V\Vert _{H^{2}(\Omega
)}^{2}$ in \eqref{3.7} by the term $\Vert V\Vert _{H^1(\Omega )}^{2}.$
This is because the $H^1(\Omega )-$norm is easier to work with
computationally than the $H^{2}(\Omega )-$norm. On the other hand, we have
not observed any instabilities probably because the number $80 \times 80$
of grid points we use is not too large and all norms in finite dimensional
spaces are equivalent.
We now identify $\{v_m(x_i, y_j): 1 \leq i, j, N_\x, 1 \leq m \leq N\}$ by the $N_\x^2 N$ dimensional vector $\mathfrak{v}$ whose $\mathfrak{i}^{\rm th}$ entry is given by
\begin{equation}
	\mathfrak{v}_{\mathfrak{i}} = v_m(x_i, y_j).
	\label{lineup}
\end{equation}
Here, $(i, j, m)$ is such that 
\begin{equation}
	\mathfrak{i} = (i - 1)N_\x N + (j - 1)N + m.
	\label{index}
\end{equation}
Then, by {approximating} all differential operators in the right hand side of \eqref{4.2} {with} their corresponding finite difference versions, we have
\begin{multline}
	J_{\epsilon}^{(0)}(V) = d_\x^2 |\mathfrak{L} \mathfrak{v}|^2
	+ d_\x|\mathfrak{D}_1 \mathfrak{v} - \mathfrak{F}|^2	
	+ d_\x|\mathfrak{D}_2 \mathfrak{v} - \mathfrak{G}|^2	
	+ \epsilon d_\x^2 |\mathfrak{v}|^2
	\\
	+ \epsilon d_\x^2 |D_x \mathfrak{v}|^2
	+ \epsilon d_\x^2 |D_y \mathfrak{v}|^2,
	\label{4.3}
\end{multline}
where the matrices $\mathfrak{L}$, $\mathfrak{D}_1$, $\mathfrak{D}_2$, $D_x$ and $D_y$ and the vectors $\mathfrak{F}$ and $\mathfrak{G}$ are described below.
The $N_\x^2N \times N_\x^2 N$ matrix $\mathfrak{L}$ is given by
\begin{enumerate}
	\item $\mathfrak{L}_{\mathfrak{i} \mathfrak{j}} = -4/d_\x^2 - s_{mm} - \Delta f(x_i, y_j)/f(x_i, y_j)$ if $\mathfrak{i} = \mathfrak{j} = (i - 1)N_\x N + (j - 1)N + m$;
	\item $\mathfrak{L}_{\mathfrak{i} \mathfrak{j}}  = 1/d_\x^2$ if  $\mathfrak{i} = (i - 1)N_\x N + (j - 1)N + m$ and $\mathfrak{j} = (i \pm 1 - 1)N_\x N +t (j \pm 1 - 1)N + m$; 
	\item $\mathfrak{L}_{\mathfrak{i} \mathfrak{j}}  = -s_{mn}$ if  $\mathfrak{i} = (i - 1)N_\x N + (j - 1)N + m$ and $\mathfrak{j} = (i  - 1)N_\x N + (j  - 1)N + n$, $n \not = m;$
	\item all other entries of $\mathfrak{L}$ are $0$;
\end{enumerate}
for all $2 \leq i, j \leq N_\x - 1,$ $1 \leq m \leq N.$
The $N_\x^2N \times N_\x^2 N$ matrix $\mathfrak{D}_1$ is given by
\begin{enumerate}
	\item $(\mathfrak{D}_1)_{\mathfrak{i} \mathfrak{j}} = 1$ if  $\mathfrak{i} = \mathfrak{j} = (i - 1)N_\x N + (j - 1)N + m$ for $i \in \{1, N_\x\}$, $1 \leq j \leq N_\x$, $1 \leq m \leq N;$
	\item $(\mathfrak{D}_1)_{\mathfrak{i} \mathfrak{j}} = 1$ if  $\mathfrak{i} = \mathfrak{j} = (i - 1)N_\x N + (j - 1)N + m$ for $1 \leq i \leq N_\x,$ $j \in \{1, N_\x\}$, $1 \leq m \leq N;$
	\item all other entries of $\mathfrak{D}_1$ are $0$.
\end{enumerate}
The $N_\x^2N \times N_\x^2 N$ matrix $\mathfrak{D}_2$ is given by
\begin{enumerate}
	\item $(\mathfrak{D}_2)_{\mathfrak{i} \mathfrak{j}} = 1/d_\x$ if  $\mathfrak{i} = \mathfrak{j} = (i - 1)N_\x N + (j - 1)N + m$ for $i \in \{1, N_\x\}$, $1 \leq j \leq N_\x$, $1 \leq m \leq N;$
	 \item $(\mathfrak{D}_2)_{\mathfrak{i} \mathfrak{j}} = -1/d_\x$ if  $\mathfrak{i}  = (i - 1)N_\x N + (j - 1)N + m$ and $\mathfrak{j}  = (i + 1 - 1)N_\x N + (j - 1)N + m$ for $i = 1$, $1 \leq j \leq N_\x$, $1 \leq m \leq N;$
	 \item $(\mathfrak{D}_2)_{\mathfrak{i} \mathfrak{j}} = -1/d_\x$ if  $\mathfrak{i}  = (i - 1)N_\x N + (j - 1)N + m$ and $\mathfrak{j}  = (i - 1 - 1)N_\x N + (j - 1)N + m$ for $i = N_\x$, $1 \leq j \leq N_\x$, $1 \leq m \leq N;$
	 \item $(\mathfrak{D}_2)_{\mathfrak{i} \mathfrak{j}} = 1/d_\x$ if  $\mathfrak{i} = \mathfrak{j} = (i - 1)N_\x N + (j - 1)N + m$ for $i 1 \leq i \leq N_\x$, $ j \in \{1, N_\x\}$, $1 \leq m \leq N;$
	 \item $(\mathfrak{D}_2)_{\mathfrak{i} \mathfrak{j}} = -1/d_\x$ if  $\mathfrak{i}  = (i - 1)N_\x N + (j - 1)N + m$ and $\mathfrak{j}  = (i - 1)N_\x N + (j + 1 - 1)N + m$ for  $1 \leq i \leq N_\x$, $j = 1$, $1 \leq m \leq N;$
	 \item $(\mathfrak{D}_2)_{\mathfrak{i} \mathfrak{j}} = -1/d_\x$ if  $\mathfrak{i}  = (i - 1)N_\x N + (j - 1)N + m$ and $\mathfrak{j}  = (i - 1)N_\x N + (j - 1 - 1)N + m$ for  $1 \leq i \leq N_\x$, $j = N_\x$, $1 \leq m \leq N;$
	 \item all other entries of $\mathfrak{D}_2$ are $0$.
\end{enumerate}
The $N_\x^2N \times N_\x^2 N$ matrix $D_x$ is given by
 \begin{enumerate}
 	\item $(D_x)_{\mathfrak{i} \mathfrak{j}} = 1/d_\x$ if  $\mathfrak{i} = \mathfrak{j} = (i - 1)N_\x N + (j - 1)N + m$ for $1 \leq j \leq N_\x-1$, $1 \leq j \leq N_\x-1$, $1 \leq m \leq N;$
	\item $(D_x)_{\mathfrak{i} \mathfrak{j}} = -1/d_\x$ if  $\mathfrak{i}  = (i - 1)N_\x N + (j - 1)N + m$ and $\mathfrak{j}  = (i + 1 - 1)N_\x N + (j - 1)N + m$ for $i = 1$, $1 \leq j \leq N_\x$, $1 \leq m \leq N;$
	\item all other entries of $D_x$ are $0$.
 \end{enumerate}
 The $N_\x^2N \times N_\x^2 N$ matrix $D_y$ is given by
 \begin{enumerate}
 	\item $(D_y)_{\mathfrak{i} \mathfrak{j}} = 1/d_\x$ if  $\mathfrak{i} = \mathfrak{j} = (i - 1)N_\x N + (j - 1)N + m$ for $1 \leq j \leq N_\x-1$, $1 \leq j \leq N_\x-1$, $1 \leq m \leq N;$
	\item $(D_y)_{\mathfrak{i} \mathfrak{j}} = -1/d_\x$ if  $\mathfrak{i}  = (i - 1)N_\x N + (j + 1- 1)N + m$ and $\mathfrak{j}  = (i + 1 - 1)N_\x N + (j - 1)N + m$ for $i = 1$, $1 \leq j \leq N_\x$, $1 \leq m \leq N;$
	\item all other entries of $D_y$ are $0$.
 \end{enumerate}
The vector $\mathfrak{F}$ is defined as
\begin{enumerate}
	\item  $\mathfrak{F}_{\mathfrak{i}}  = F_m(x_i, y_j)$ if  $\mathfrak{i} =   (i - 1)N_\x N + (j - 1)N + m$ for $i \in \{1, N_\x\}$, $1 \leq j \leq N_\x$, $1 \leq m \leq N;$
	\item all other entries of $\mathfrak{F}$ are $0$.
\end{enumerate}
The vector $\mathfrak{G}$ is defined as
\begin{enumerate}
	\item  $\mathfrak{G}_{\mathfrak{i}}  = G_m(x_i, y_j)$ if  $\mathfrak{i}  = (i - 1)N_\x N + (j - 1)N + m$ for $i \in \{1, N_\x\}$, $1 \leq j \leq N_\x$, $1 \leq m \leq N;$
	\item all other entries of $\mathfrak{G}$ are $0$.
\end{enumerate}

Since $\epsilon$ is small (in our computational program $\epsilon = 10^{-9}$), to find the minimizer $V^{(0)}$ of the finite difference version of $J_{\epsilon}$ in \eqref{4.3}, we solve the linear system
\[
	\left(\left[
		\begin{array}{c}
			\mathfrak{L}\\
			\mathfrak{D}_1\\
			\mathfrak{D}_2\\			
		\end{array}
	\right]^{\rm {T}}
	\left[
		\begin{array}{c}
			\mathfrak{L}\\
			\mathfrak{D}_1\\
			\mathfrak{D}_2\\			
		\end{array}
	\right] 
	+ \epsilon (\Id + D_x^{\rm T} D_x + Dy^{\rm T} D_y) \right) \mathfrak{v} 
	= \left[
		\begin{array}{c}
			\mathfrak{L}\\
			\mathfrak{D}_1\\
			\mathfrak{D}_2\\			
		\end{array}
	\right]^{\rm {T}}
	\left[
		\begin{array}{c}
			{\bf 0}\\
			\mathfrak{F}\\
			\mathfrak{G}\\			
		\end{array}
	\right].
\]
Having $\mathfrak{v}$ in hand, we can compute $V^{(0)} = (v_1^{0}, \dots, v_N^{(0)})^{\rm T}$ using \eqref{lineup} and \eqref{index}.

\noindent{\bf  Step \ref{reg 2}.} The implementation for this step is similar to that for Step \ref{reg 1},
{namely}, we minimize
\begin{multline}
	J_{\epsilon}^{(p + 1)}(V) = 
	\sum_{m = 1}^N \int_\Omega 
	\Big|
		\Delta v_m(\x) - \sum_{n = 1}^N s_{mn} v_{n}(\x) - \frac{\Delta f(\x)}{f(\x)} v_m(\x) 
		\\
	+ \sum_{n = 1}^N \frac{\Psi_n(0)}{f(\x)}v_n(\x) v_m^{(p)}(\x)
	\Big|
	\\
	+ \int_{\partial \Omega} |\partial_{\nu}v_m(\x) - G_m(\x)|^2d\sigma(\x)
	+ \epsilon \int_{\Omega} (|v_m(\x)|^2 + |\nabla v_m(\x)|^2) d\x
	\Big].
\end{multline}
To this end, we identify the vector valued function $V$ by the vector $\mathfrak{v}$ as in \eqref{lineup} and \eqref{index} and then solve the linear system 
\[
	\left(\left[
		\begin{array}{c}
			\mathcal{L}\\
			\mathfrak{D}_1\\
			\mathfrak{D}_2\\			
		\end{array}
	\right]^{\rm T}
	\left[
		\begin{array}{c}
			\mathcal{L}\\
			\mathfrak{D}_1\\
			\mathfrak{D}_2\\			
		\end{array}
	\right] 
	+ \epsilon (\Id + D_x^{\rm T} D_x + Dy^{\rm T} D_y) \right) \mathfrak{v} 
	= \left[
		\begin{array}{c}
			\mathcal{L}\\
			\mathfrak{D}_1\\
			\mathfrak{D}_2\\			
		\end{array}
	\right]^{\rm T}
	\left[
		\begin{array}{c}
			{\bf 0}\\
			\mathfrak{F}\\
			\mathfrak{G}\\			
		\end{array}
	\right].
\]
Here, the matrices $\mathfrak{D}_1$, $\mathfrak{D}_2$, $D_x$ and $D_y$ and the vectors $\mathfrak{F}$ and $\mathfrak{G}$ are defined in the implementation section for Step \ref{reg 1}.
The $N_\x^2N \times N_\x^2 N$ matrix $\mathcal{L}$ is given by
\begin{enumerate}
	\item $\mathcal{L}_{\mathfrak{i} \mathfrak{j}} = -4/d_\x^2 - s_{mm} - \Delta f(x_i, y_j)/f(x_i, y_j) + v_m^{(p)}(x_i, y_j) \Psi_m(0)/f(x_i, y_j) $ if $\mathfrak{i} = \mathfrak{j} = (i - 1)N_\x N + (j - 1)N + m$;
	\item $\mathcal{L}_{\mathfrak{i} \mathfrak{j}}  = 1/d_\x^2$ if $\mathfrak{i} = (i - 1)N_\x N + (j - 1)N + m$ and $\mathfrak{j} = (i \pm 1 - 1)N_\x N + (j \pm 1 - 1)N + m$; 
	\item $\mathcal{L}_{\mathfrak{i} \mathfrak{j}}  = -s_{mn} +  v_n^{(p)}(x_i, y_j) \Psi_n(0)/f(x_i, y_j)$ if  $\mathfrak{i} = (i - 1)N_\x N + (j - 1)N + m$ and $\mathfrak{j} = (i  - 1)N_\x N + (j  - 1)N + n $, $n \not = m;$
	\item all other entries of $\mathcal{L}$ are $0$;
\end{enumerate}
for all $2 \leq i, j \leq N_\x - 1,$ $1 \leq m \leq N.$
Having $\mathfrak{v}$ in hand, we can compute $V^{(p+1)} = (v_1^{(p + 1)}, \dots, v_N^{(p + 1)})^T$ using \eqref{lineup} and \eqref{index}.

\noindent {\bf Steps \ref{st3}, \ref{st5} and \ref{st6}.}  The implementation of these steps is straight forward. 

{
\begin{remark}
	All matrices above are of the large size $N_\x^2 N \times N_\x^2 N$, which might cause some inefficiency in computations. 
	However, since most of their entries are zeros, we can treat those matrices as sparse ones to overcome this difficulty.
	The linear algebra package for sparse matrices are already built in Matlab.
\end{remark}
}

In the next section, we show some numerical results.

\section{Numerical examples} \label{sec numerical example}

The numerical results presented below are computed from the knowledge of $F_m(\x)$ and $G_m(\x)$, $m \in \{1, \dots, N\}$, on $\partial \Omega \times [0, 0.3]$ including 10\% of noise, where $F$ and $G$ are the boundary data in Remark \ref{rem noise}.
The number of truncation $N$ is $25$.
The regularization parameter is $\epsilon = 10^{-9}.$
The computational program is implemented by the finite difference method.

\begin{enumerate}
\item {\it Test 1.} The true function $c_{\rm true}$ has a smooth inclusion
\[
	c_{\rm true}{(x, y)} = \left\{
		\begin{array}{ll}
			0 & x^2 + (y + 0.3)^2 \geq 0.35^2\\
			20e^{\frac{x^2 + (y + 0.3)^2}{x^2 + (y + 0.3)^2 - 0.23^2}} &x^2 + (y + 0.3)^2 \geq 0.35^2.
		\end{array}
	\right.
\]

\begin{figure}
	\subfloat[The function $c_{\rm true}$]{
		\includegraphics[width=.3\textwidth]{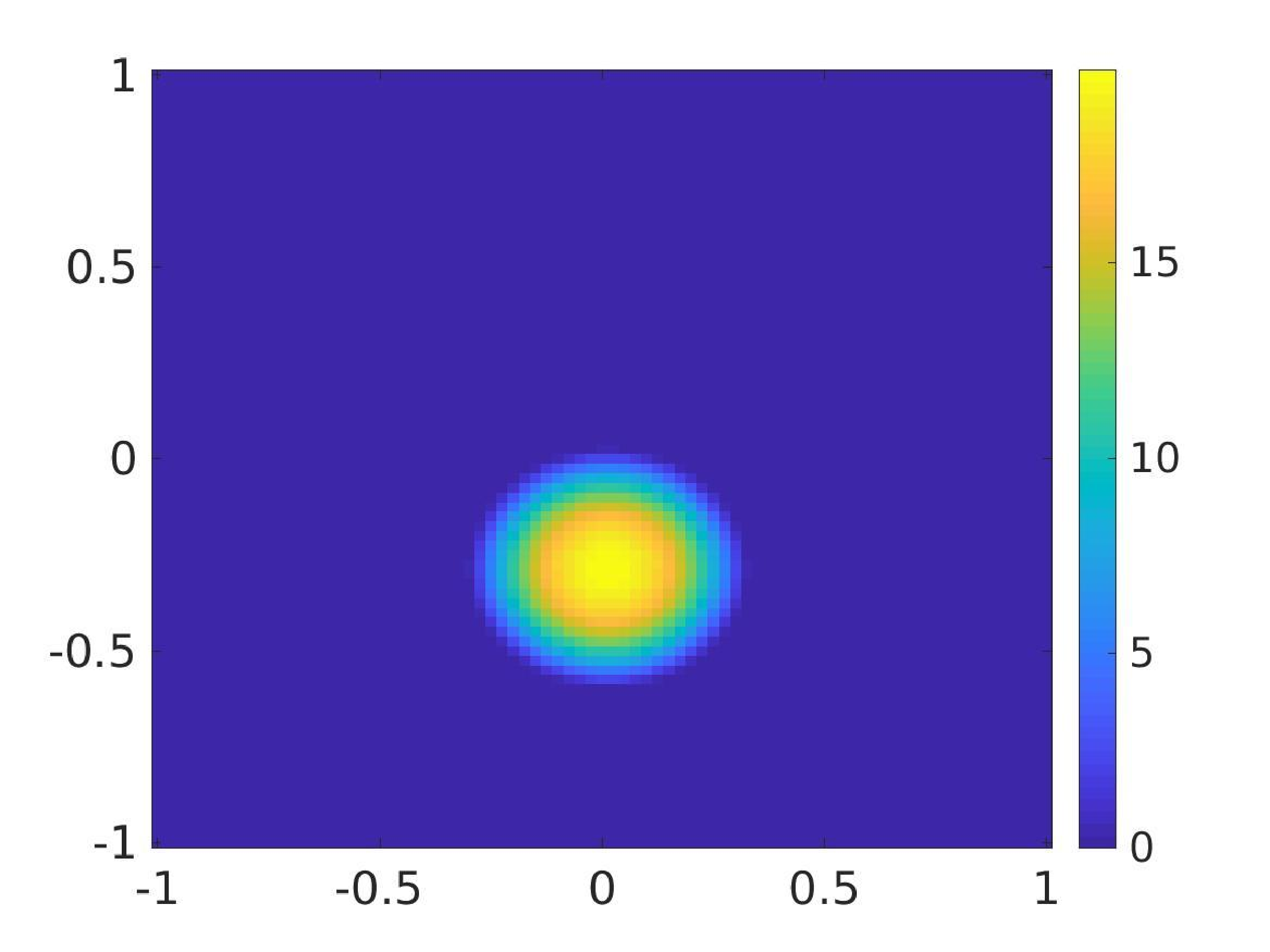}
	}  	
	\subfloat[\label{1b}The function $c^{(0)}$]{
		\includegraphics[width=.3\textwidth]{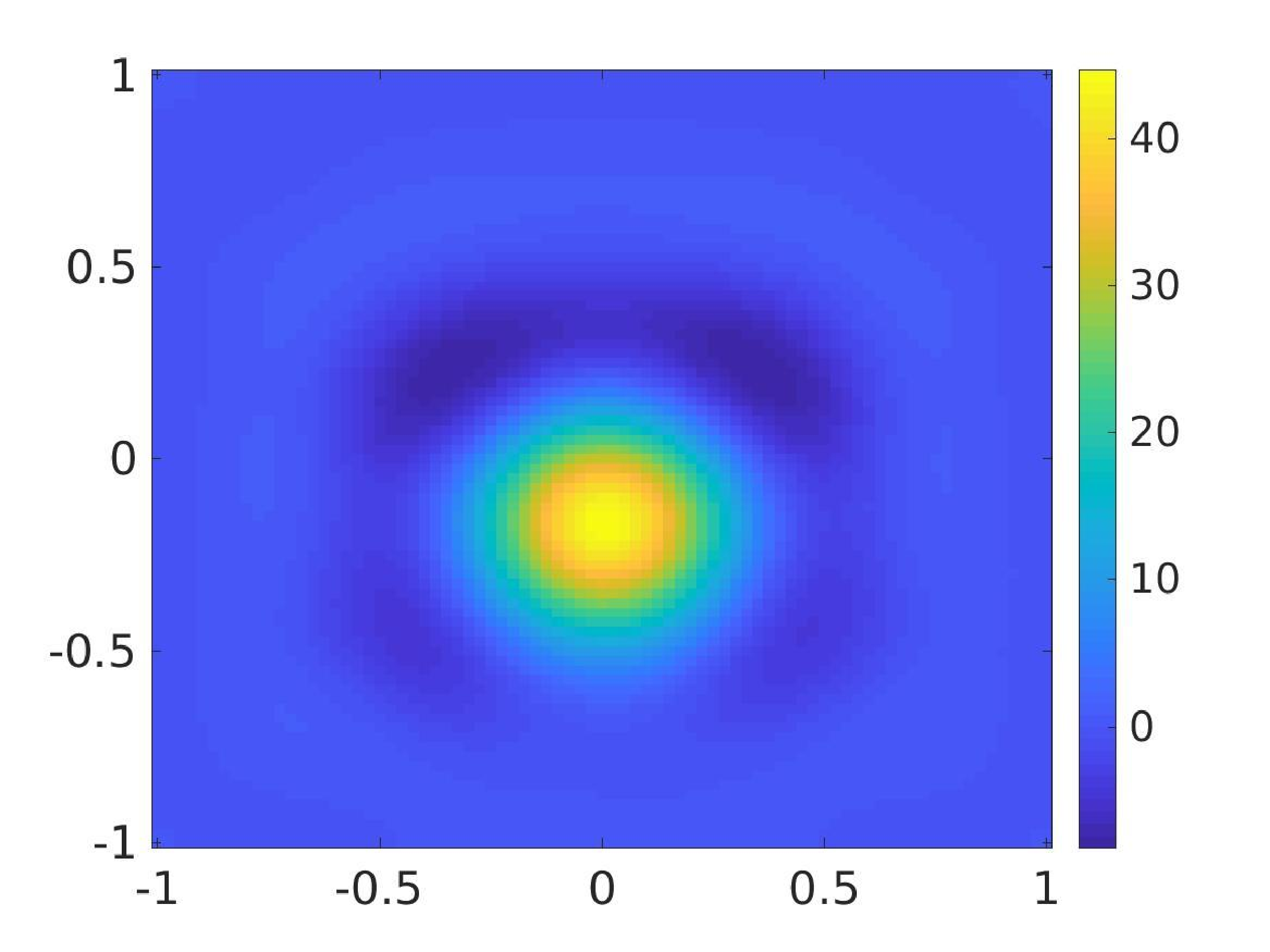}
	} 
	\subfloat[\label{1c}The function $c^{(1)}$]{
		\includegraphics[width=.3\textwidth]{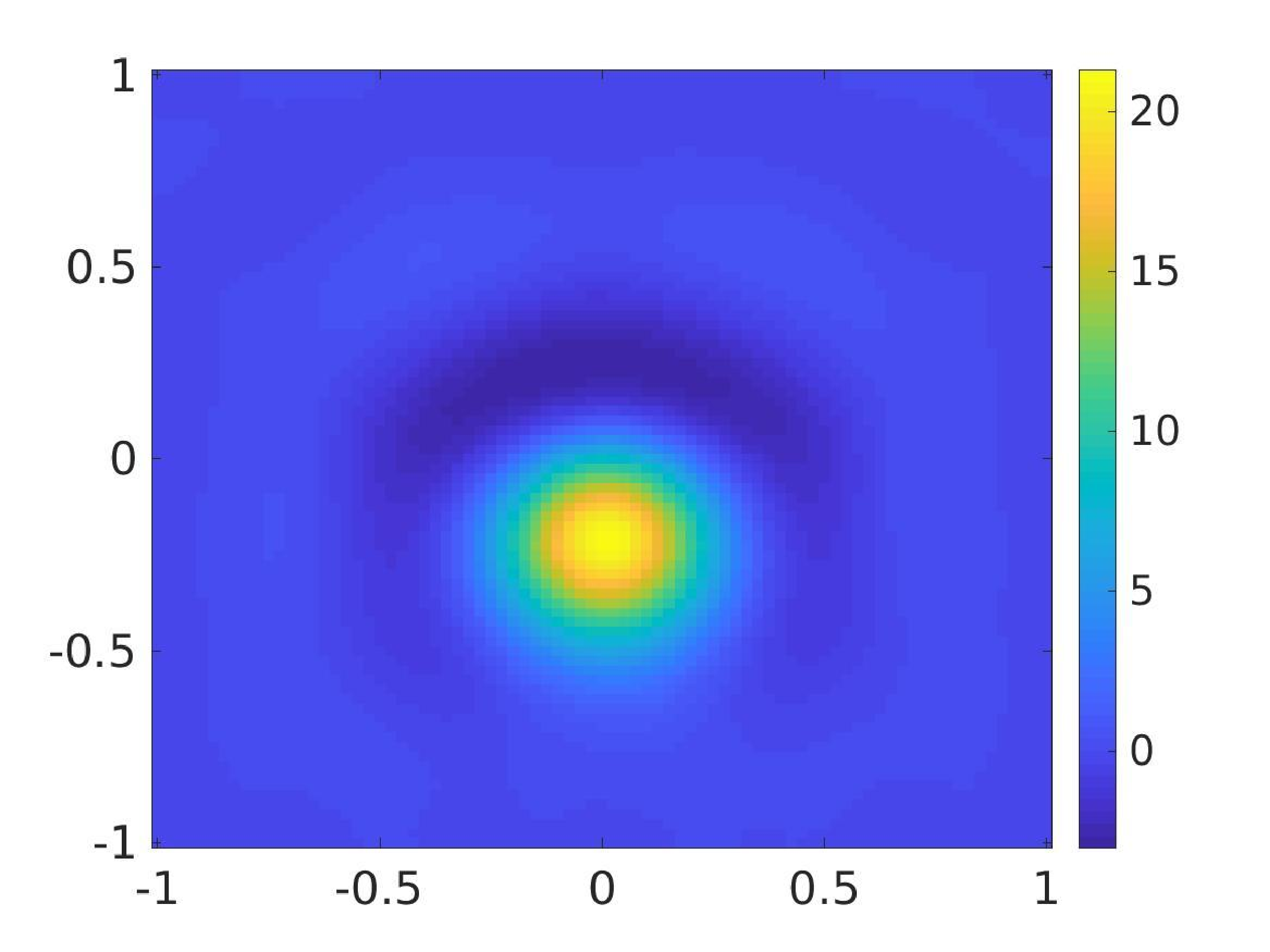}
		}	
		
	\subfloat[The function $c^{(2)}$]{
		\includegraphics[width=.3\textwidth]{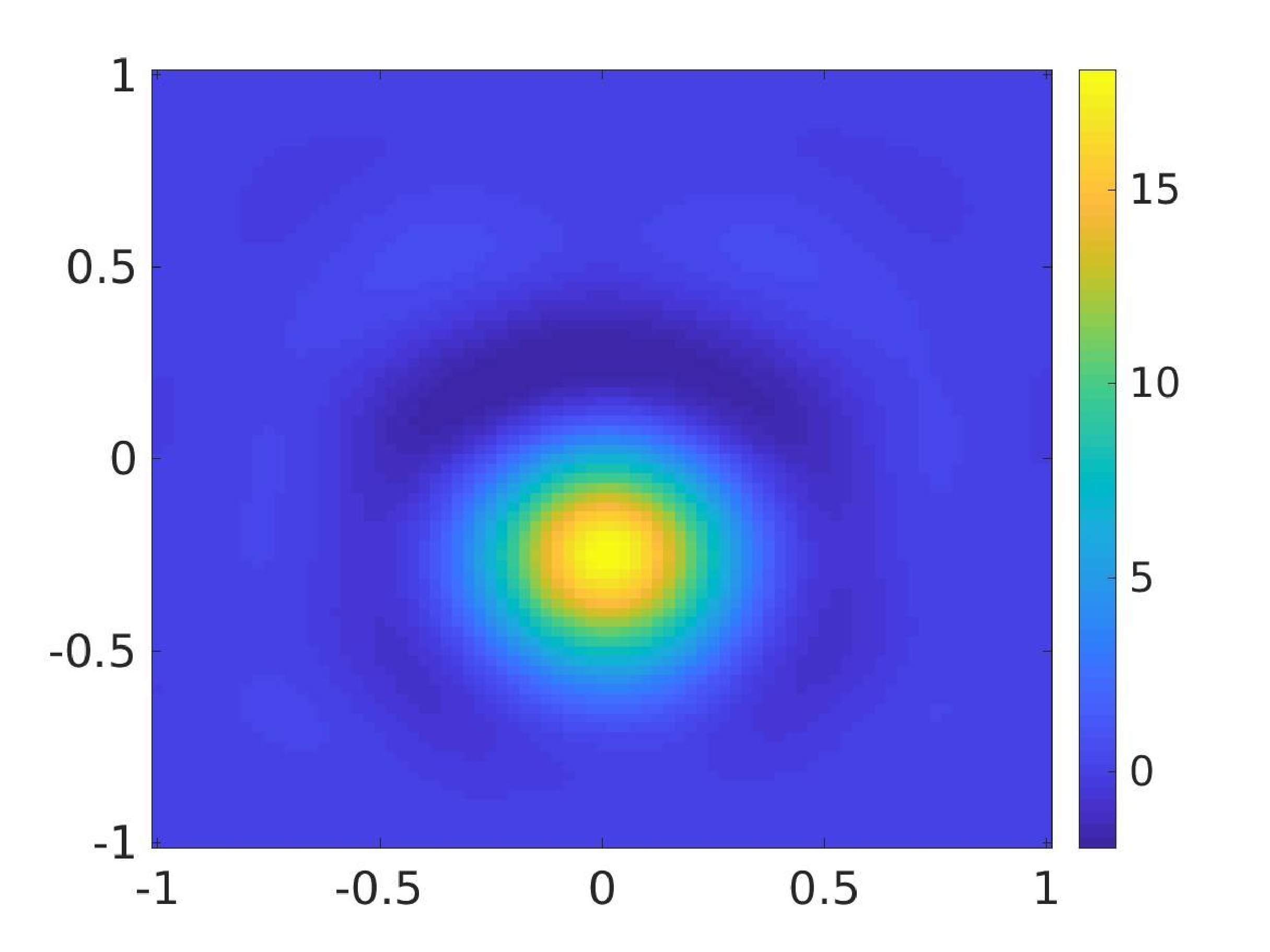}
	}			
	\subfloat[\label{1d}The function $c^{(10)}$]{
		\includegraphics[width=.3\textwidth]{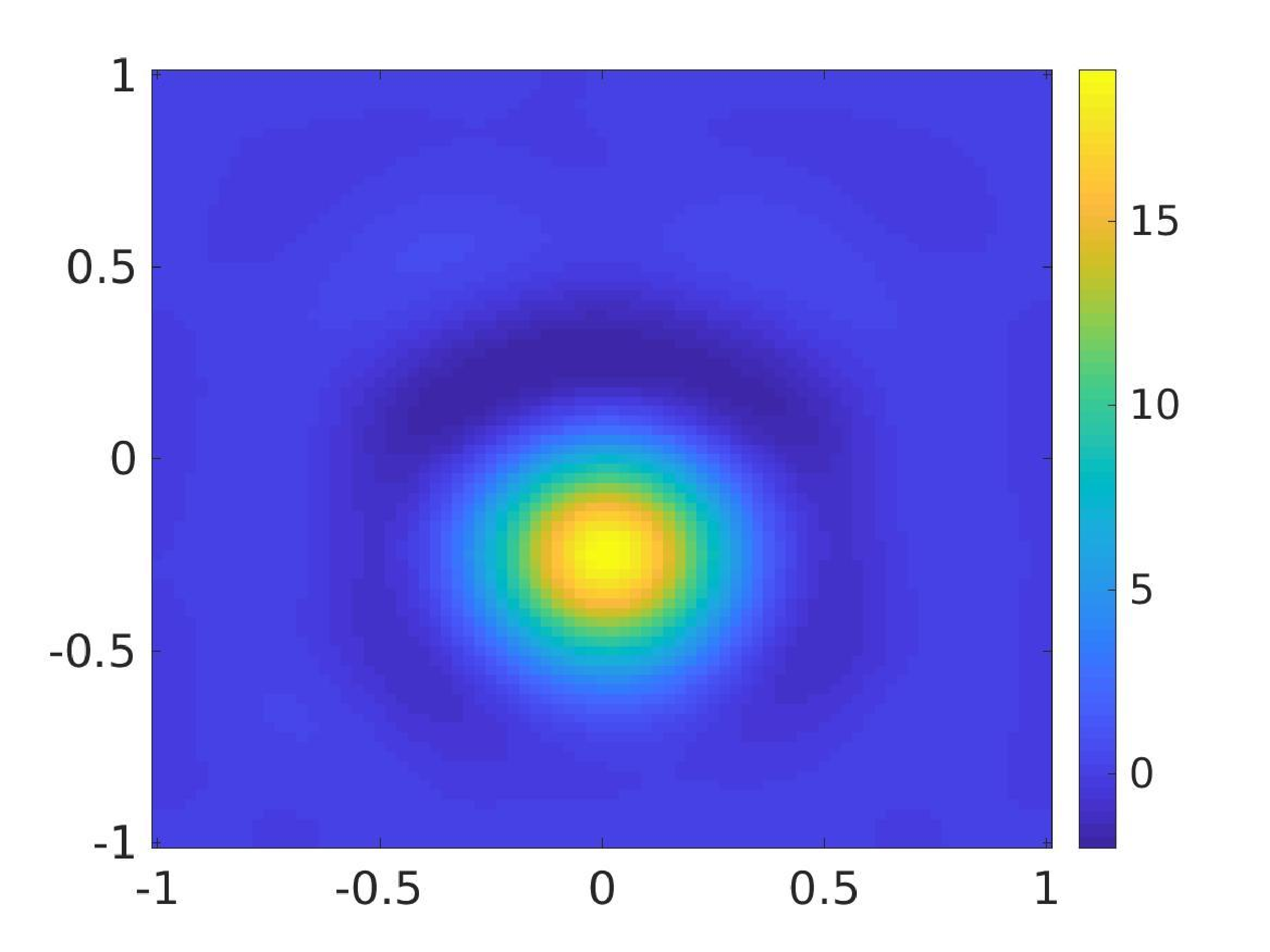}
	}		
	\subfloat[\label{1e}The error function ${E}$]{
		\includegraphics[width=.3\textwidth]{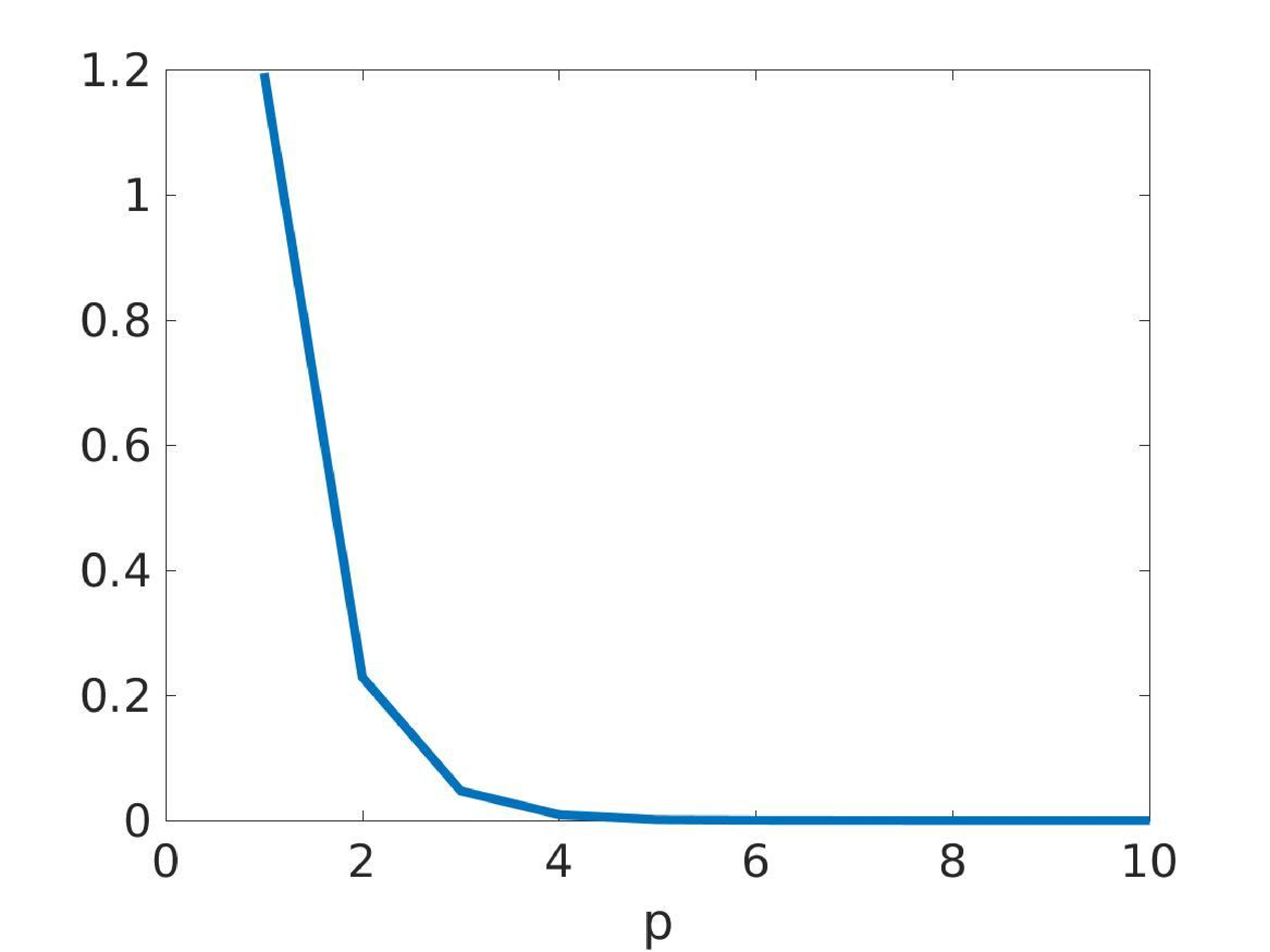}
	}

	\caption{\label{fig 1}Test 1. The true coefficient  and computed coefficient $c$. We observe from Figure \ref{1b} that $c^{(0)}$, computed by Step \ref{reg 1} of Algorithm \ref{alg}, has good ``image" of the true inclusion. 
	It is evident from the graph of the error function ${E}$, see Figure \ref{1e}, that the sequence $\{c^{(p)}\}_{{p \geq 1}}$ converges fast.
	}
\end{figure}

The numerical results for this case {are} displayed in Figure \ref{fig 1}. 
One can observe in Figures \ref{1b}--\ref{1d} that the circular shape and location of the inclusion can be succesfully detected.
The true maximal value of the function $c_{\rm true}$ is 20.
The reconstructed maximal value of the function $c_{\rm comp} = c^{(10)}$ is $19.07$. 

\item {\it Test 2.} We test the case when the function $c_{\rm true}$ is a step function with two rectangular inclusions. 
This example is interesting since $c_{\rm true}$ is not smooth and 
the gap at the boundaries  of the inclusions is high.
The function $c_{\rm true}$ is given by
\[
	c_{\rm true}{(x, y)} = \left\{
		\begin{array}{ll}
			10 &|x| < 0.8 \mbox{ and } |y \pm 0.4| < 0.15,\\
			0 &\mbox{otherwise}.
		\end{array}
	\right.
\]

\begin{figure}
	\subfloat[The function $c_{\rm true}$]{
		\includegraphics[width=.3\textwidth]{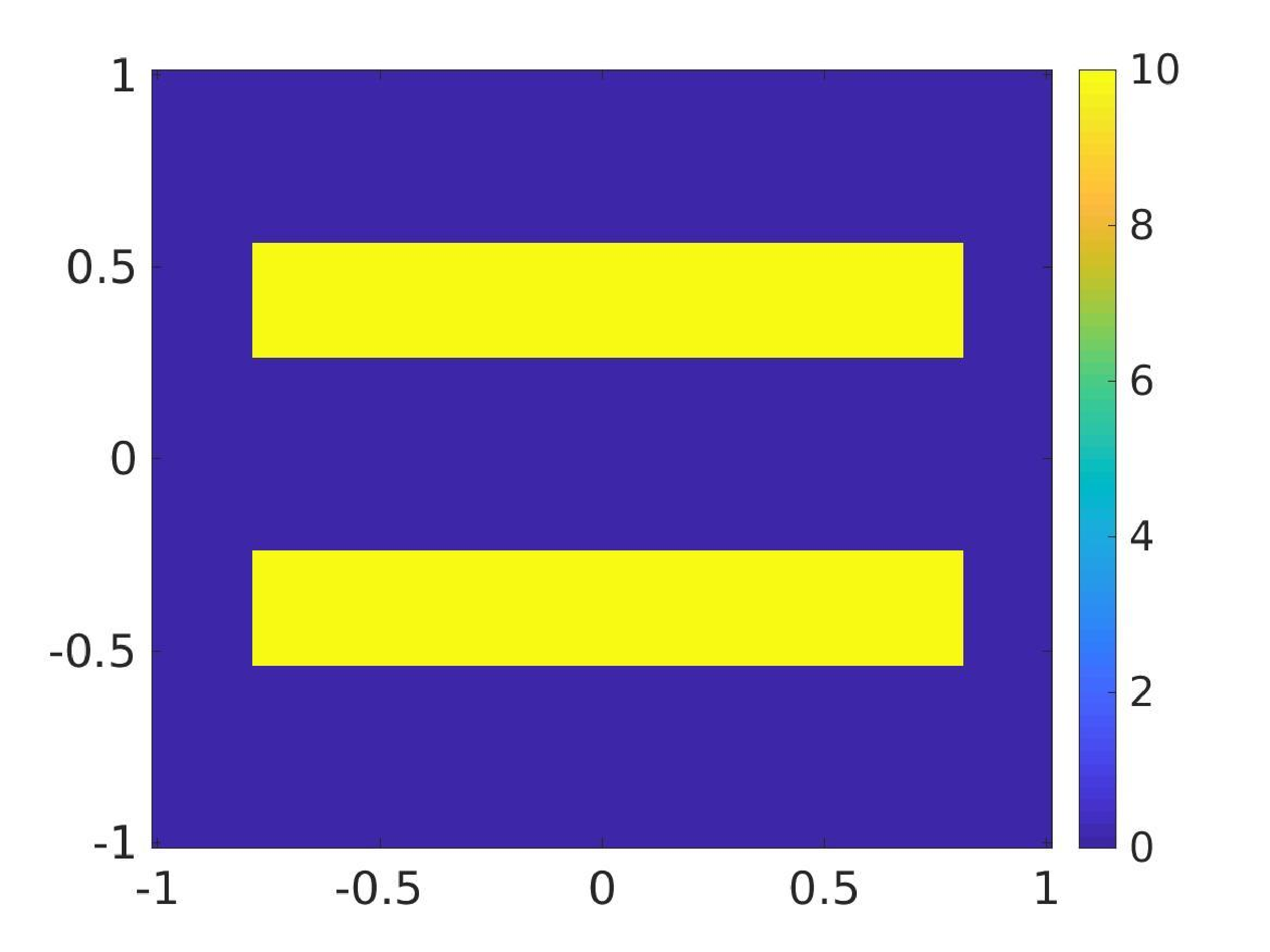}
	}  	
	\subfloat[\label{2b}The function $c^{(0)}$]{
		\includegraphics[width=.3\textwidth]{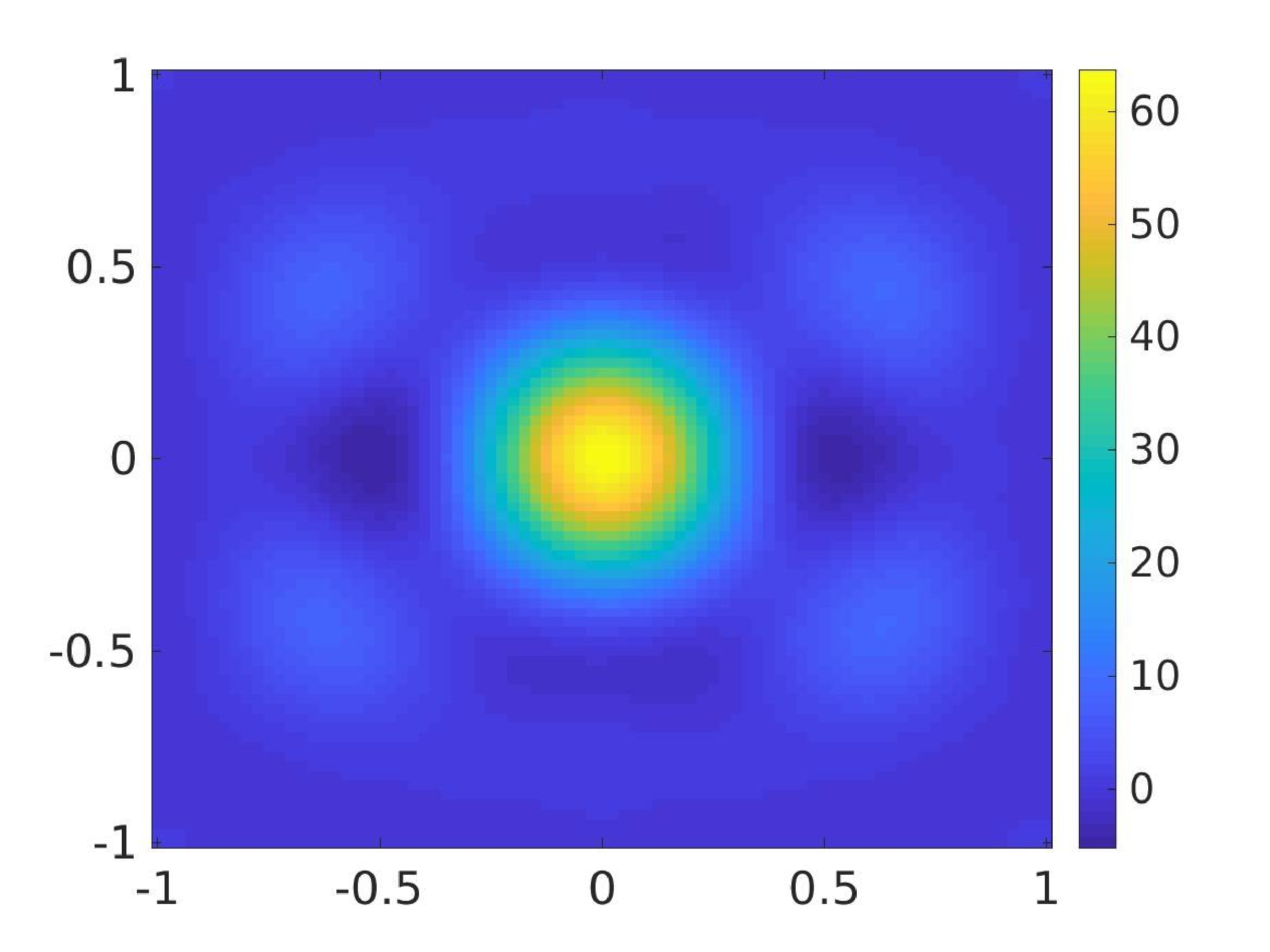}
	} 
	\subfloat[\label{2c}The function $c^{(1)}$]{
		\includegraphics[width=.3\textwidth]{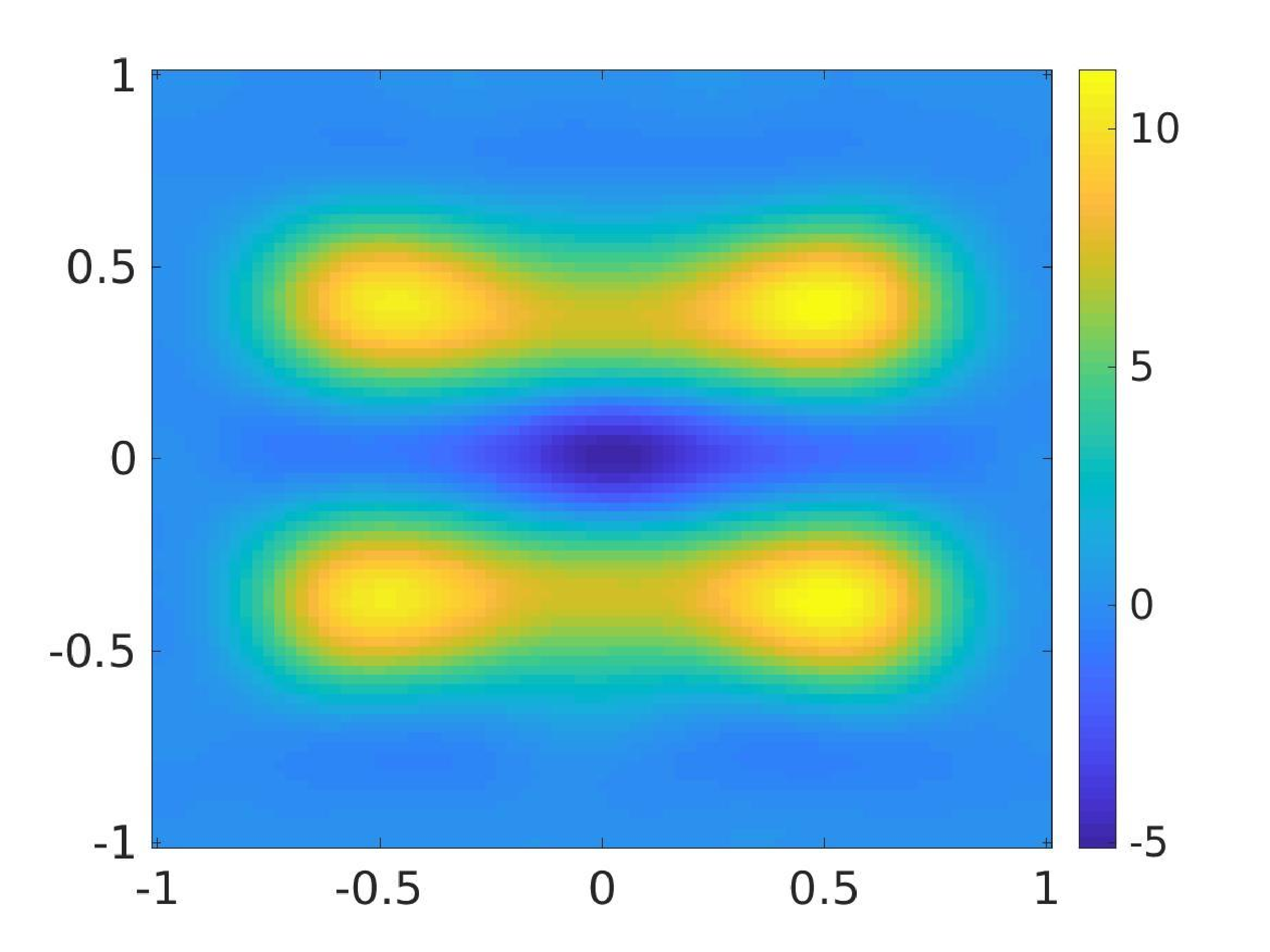}
		}

	\subfloat[The function $c^{(2)}$]{
		\includegraphics[width=.3\textwidth]{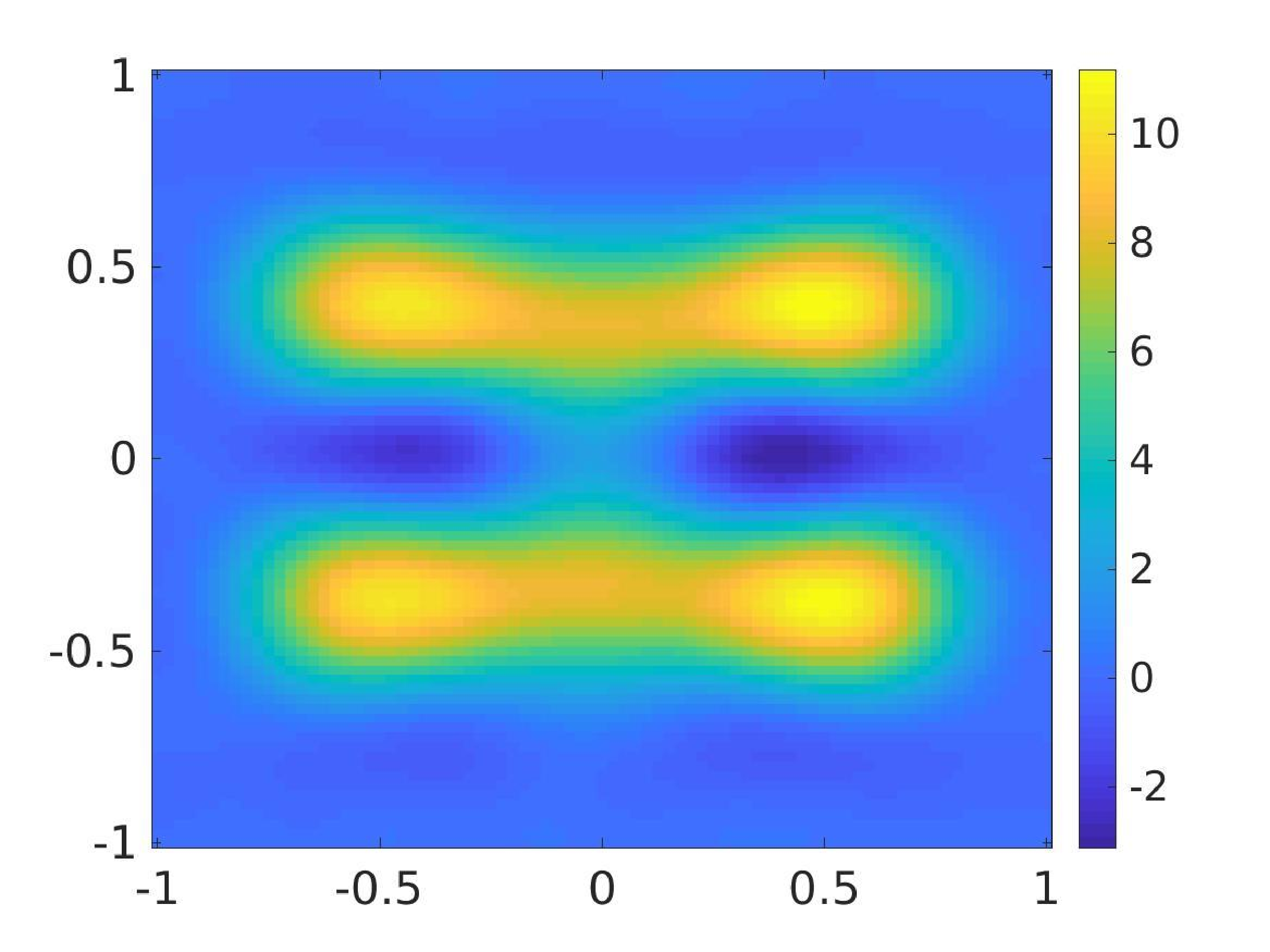}
	}	
	\subfloat[\label{2d}The function $c^{(10)}$]{
		\includegraphics[width=.3\textwidth]{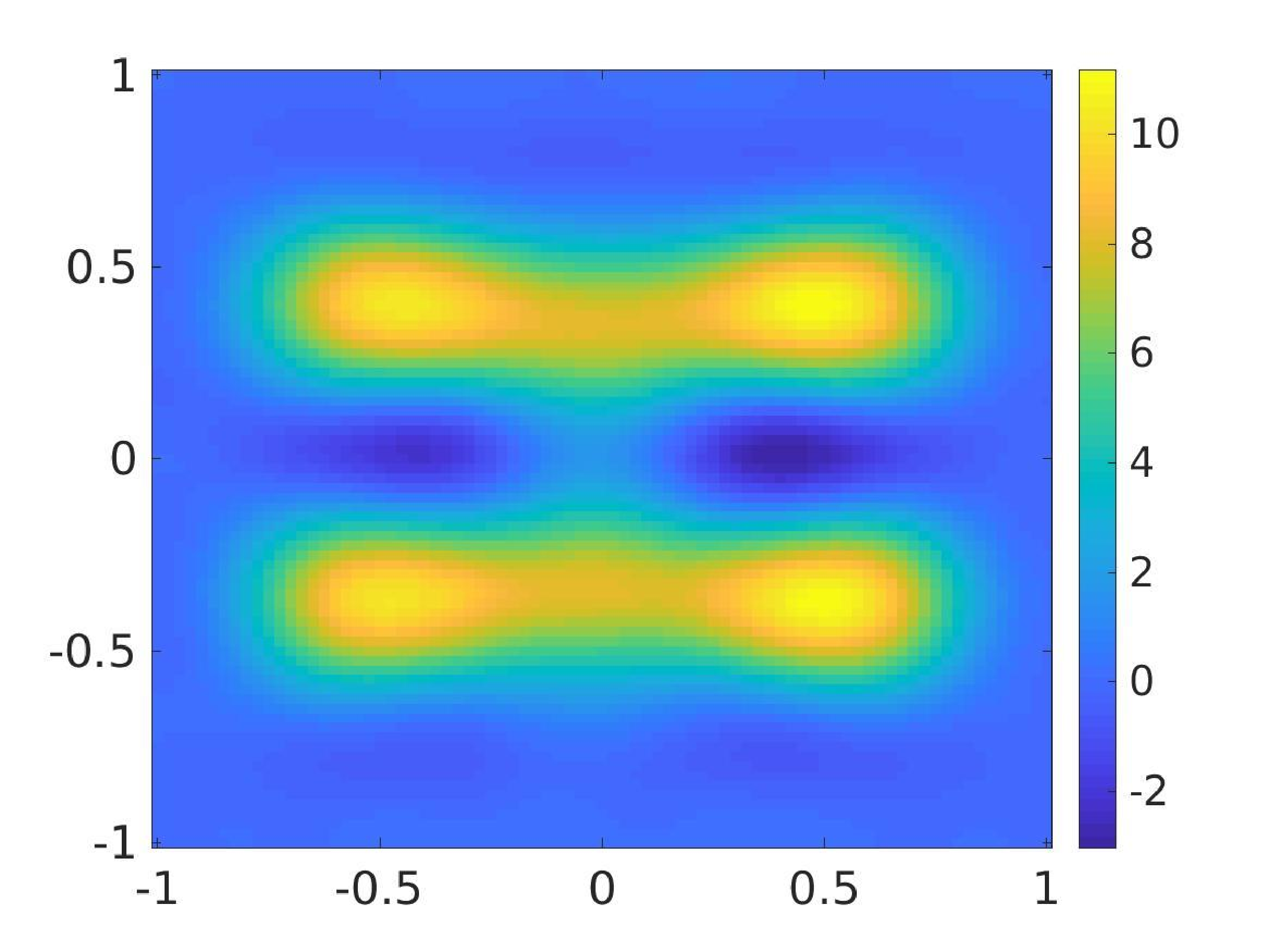}
	}		
	\subfloat[\label{2e}The error function ${E}$]{
		\includegraphics[width=.3\textwidth]{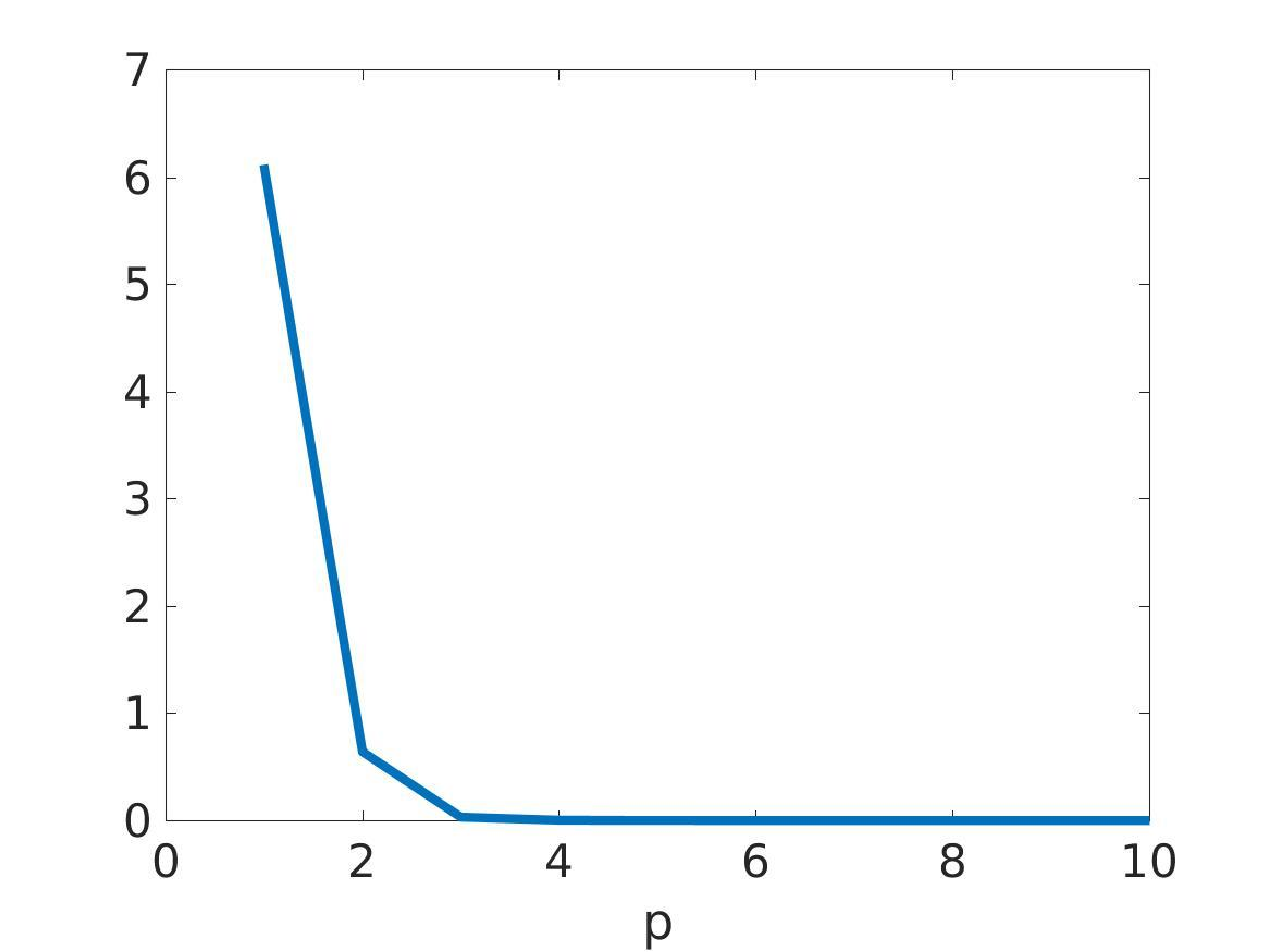}
	}

	\caption{\label{fig 2}Test 2. The true coefficient  and computed coefficient $c$. 
	In this test, although the reconstructed coefficient $c^{(0)}$, see Figure \ref{2b}, is poor, the reconstructed coefficient $c^{(1)}$ meets the expectation.
	It is evident from the graph of the error function ${E}$, see Figure \ref{2e}, that the sequence $\{c^{(p)}\}_{{p \geq 1}}$ converges fast.
	}
\end{figure}

The numerical results for this case {are} displayed in Figure \ref{fig 2}. 
One can observe in Figures \ref{2c}--\ref{2d} that the reconstructed rectangular shape and location of the inclusion are satisfactory.
The true maximal value of the function $c_{\rm true}$ is 10.
The reconstructed maximal value of the function $c_{\rm comp}$ is 10.98. 
Similarly to the previous test, it is evident from Figure \ref{2e}  that our method converges fast.

\item {\it Test 3.} 
We test the case of two circular inclusions. 
In this case, the function $c_{\rm comp}$ is a step function with high gap at the boundary of the inclusions.
The function $c_{\rm true}$ is given by
\[
	c_{\rm true}{(x, y)} = \left\{
		\begin{array}{ll}
			5 &x^2 + (y + 0.5)^2 < 0.23^2,\\
			8 &x^2 + (y - 0.5)^2 < 0.23^2,\\
			0 &\mbox{otherwise}.
		\end{array}
	\right.
\]

\begin{figure}[h!]
	\subfloat[The function $c_{\rm true}$]{
		\includegraphics[width=.3\textwidth]{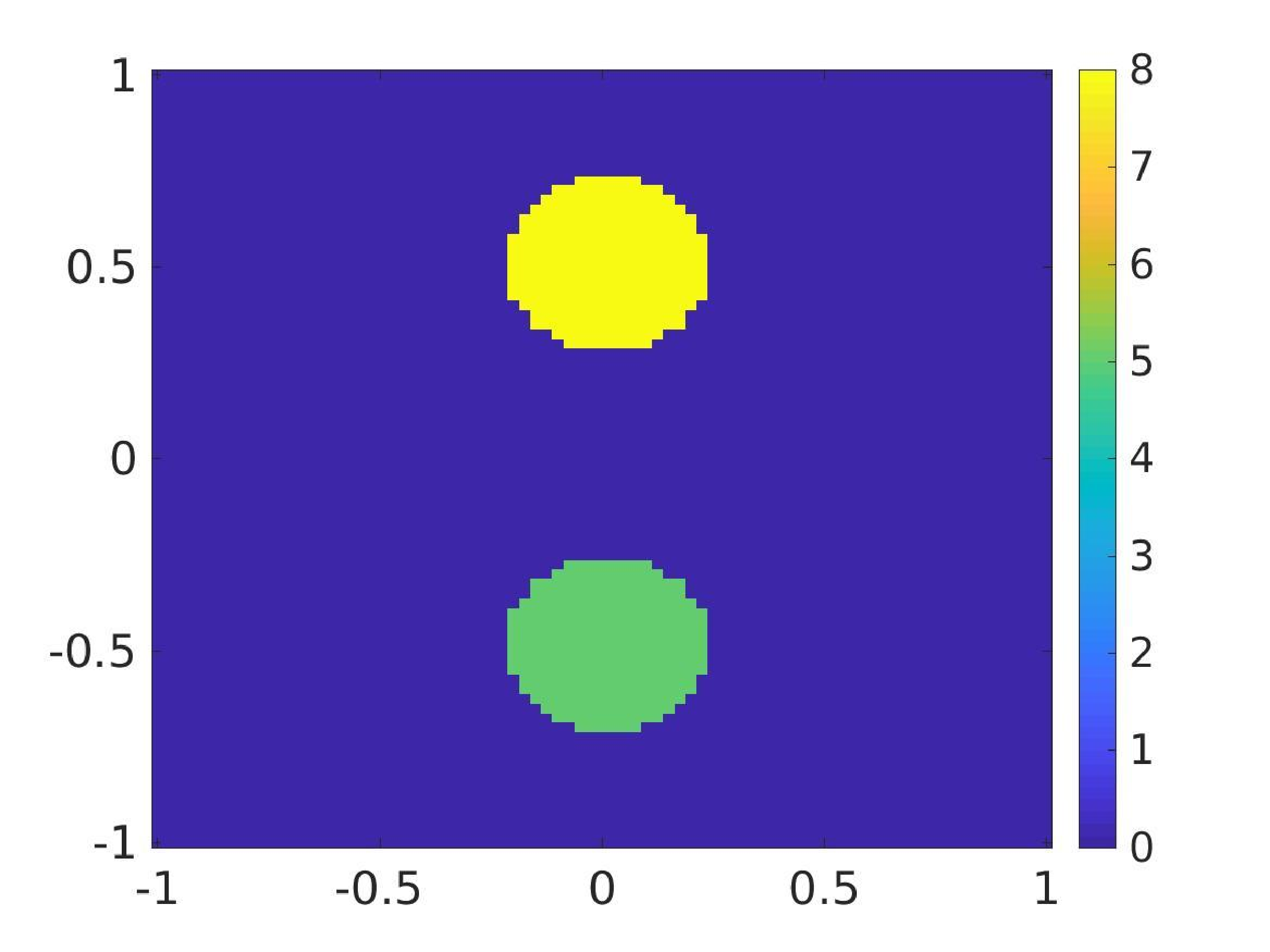}
	}  	
	\subfloat[\label{3b}The function $c^{(0)}$]{
		\includegraphics[width=.3\textwidth]{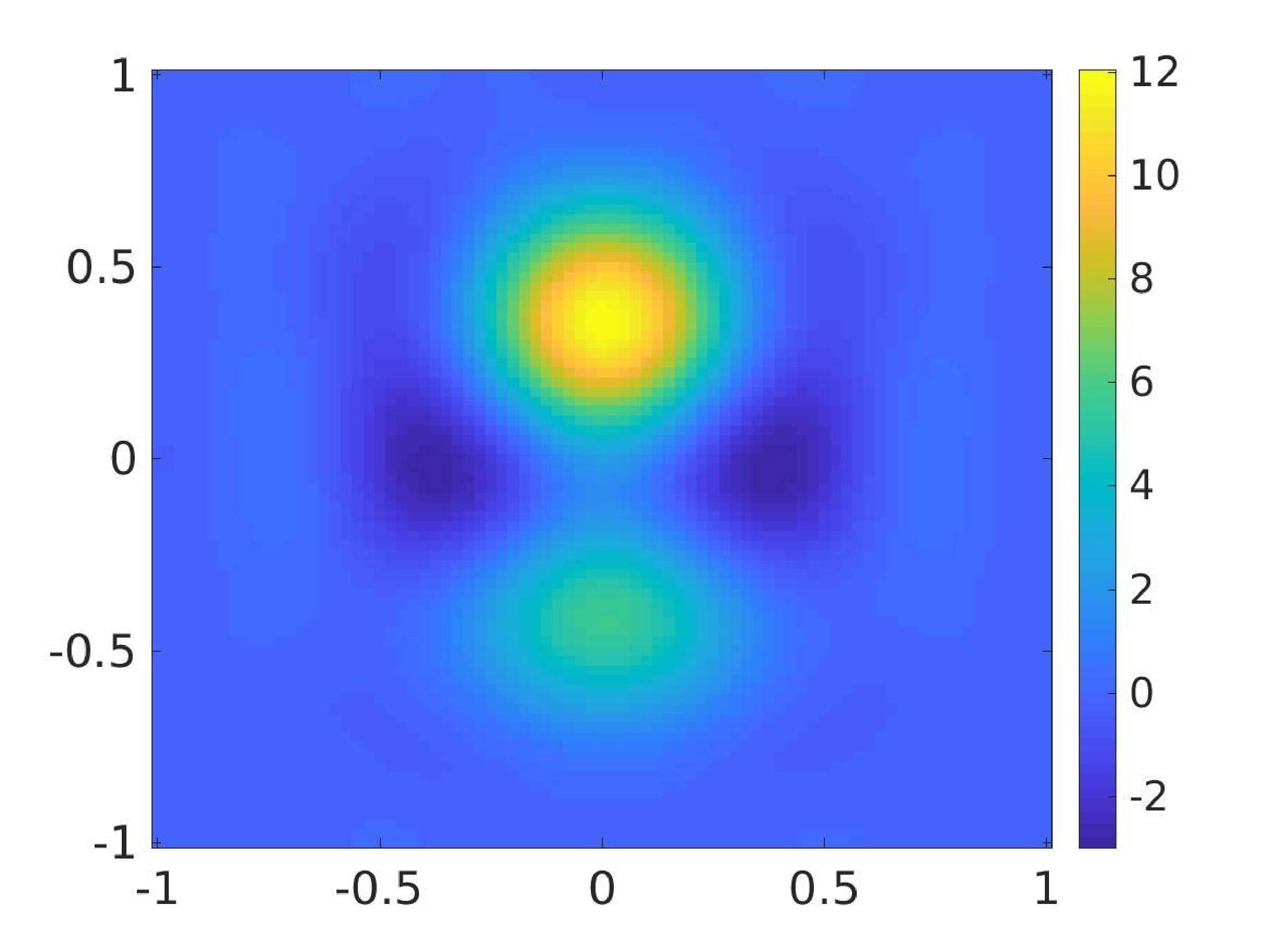}
	} 
	\subfloat[\label{3c}The function $c^{(1)}$]{
		\includegraphics[width=.3\textwidth]{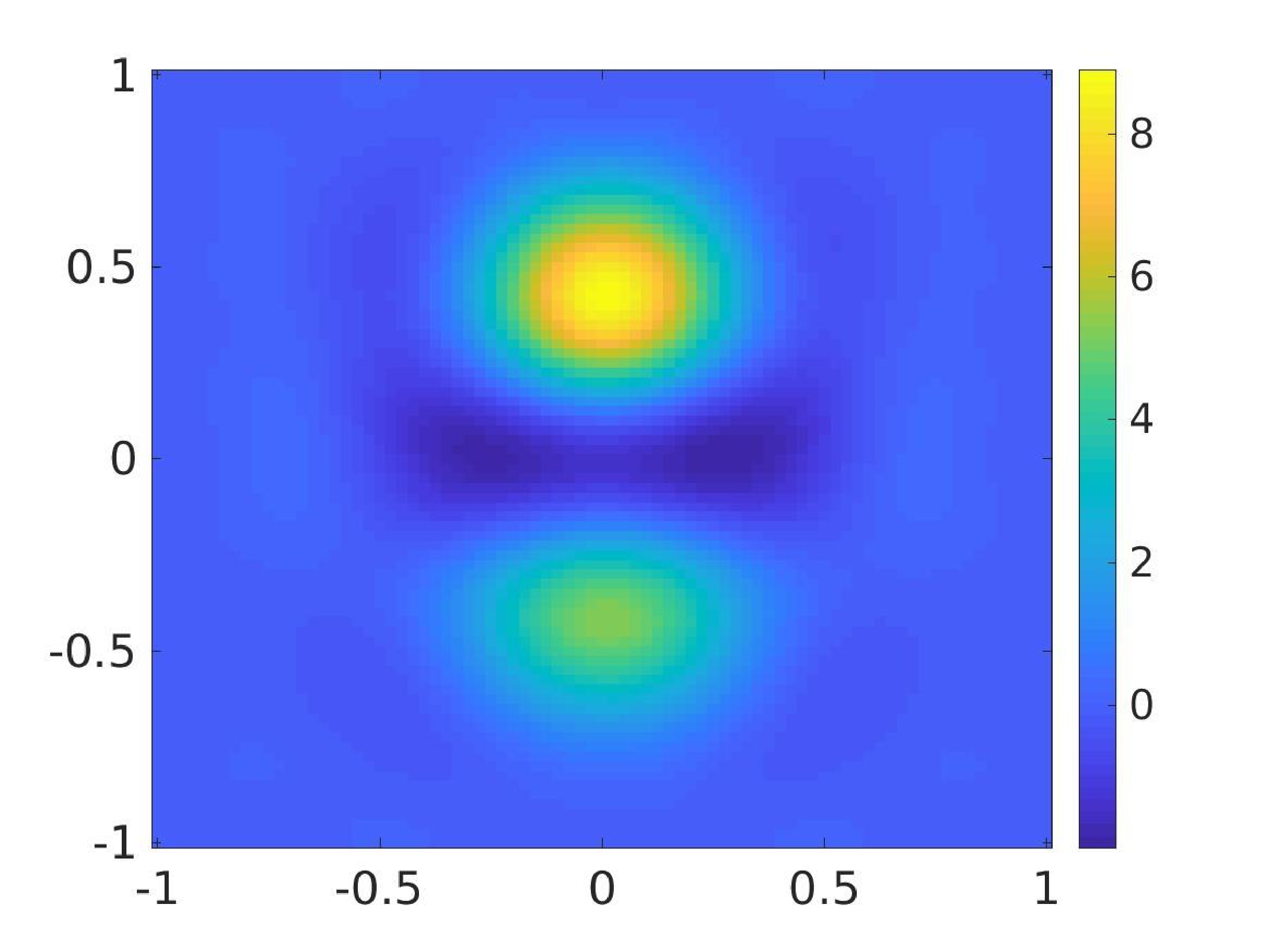}
		}

	\subfloat[The function $c^{(2)}$]{
		\includegraphics[width=.3\textwidth]{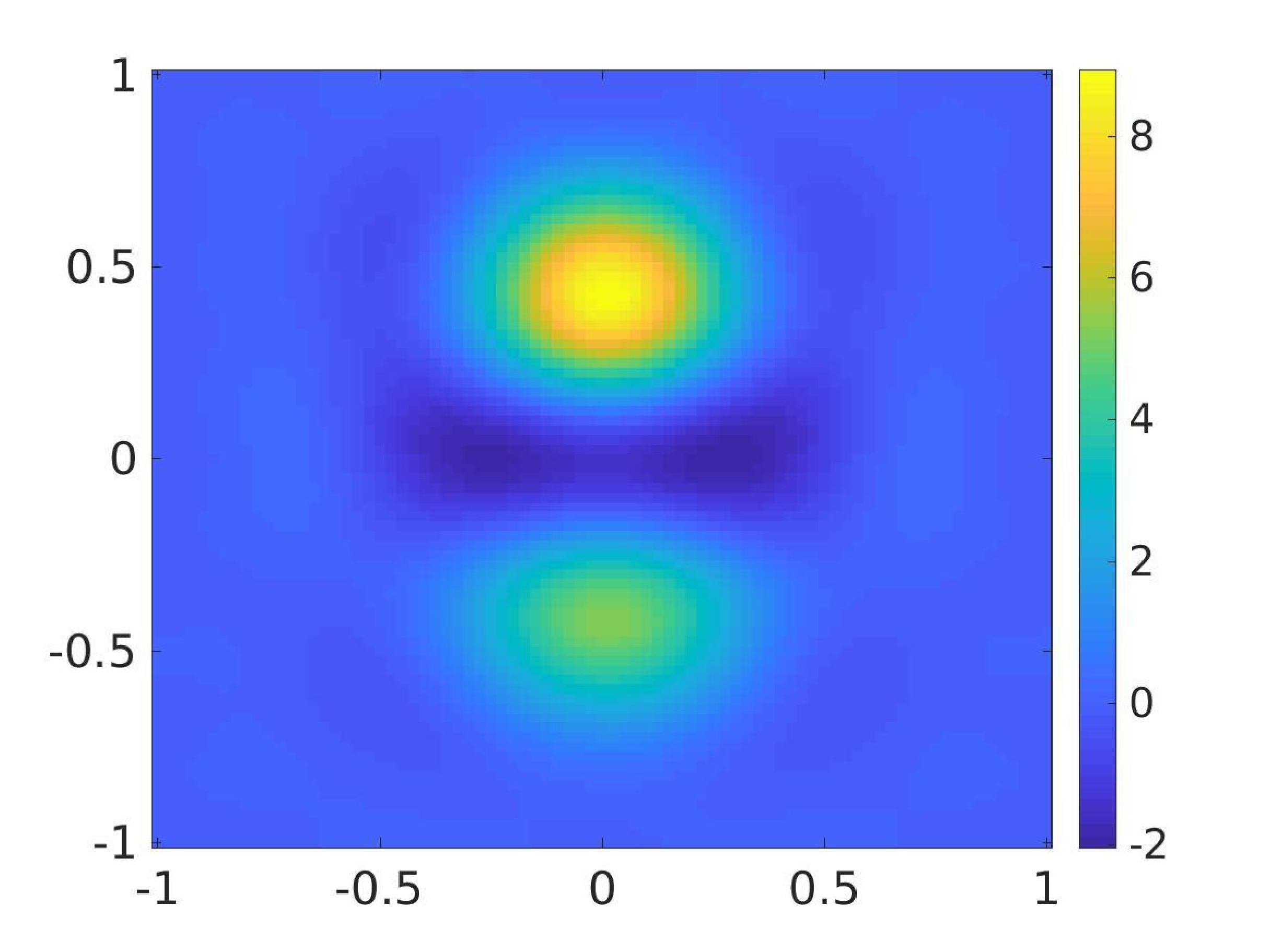}
	}	
	\subfloat[\label{3d}The function $c^{(10)}$]{
		\includegraphics[width=.3\textwidth]{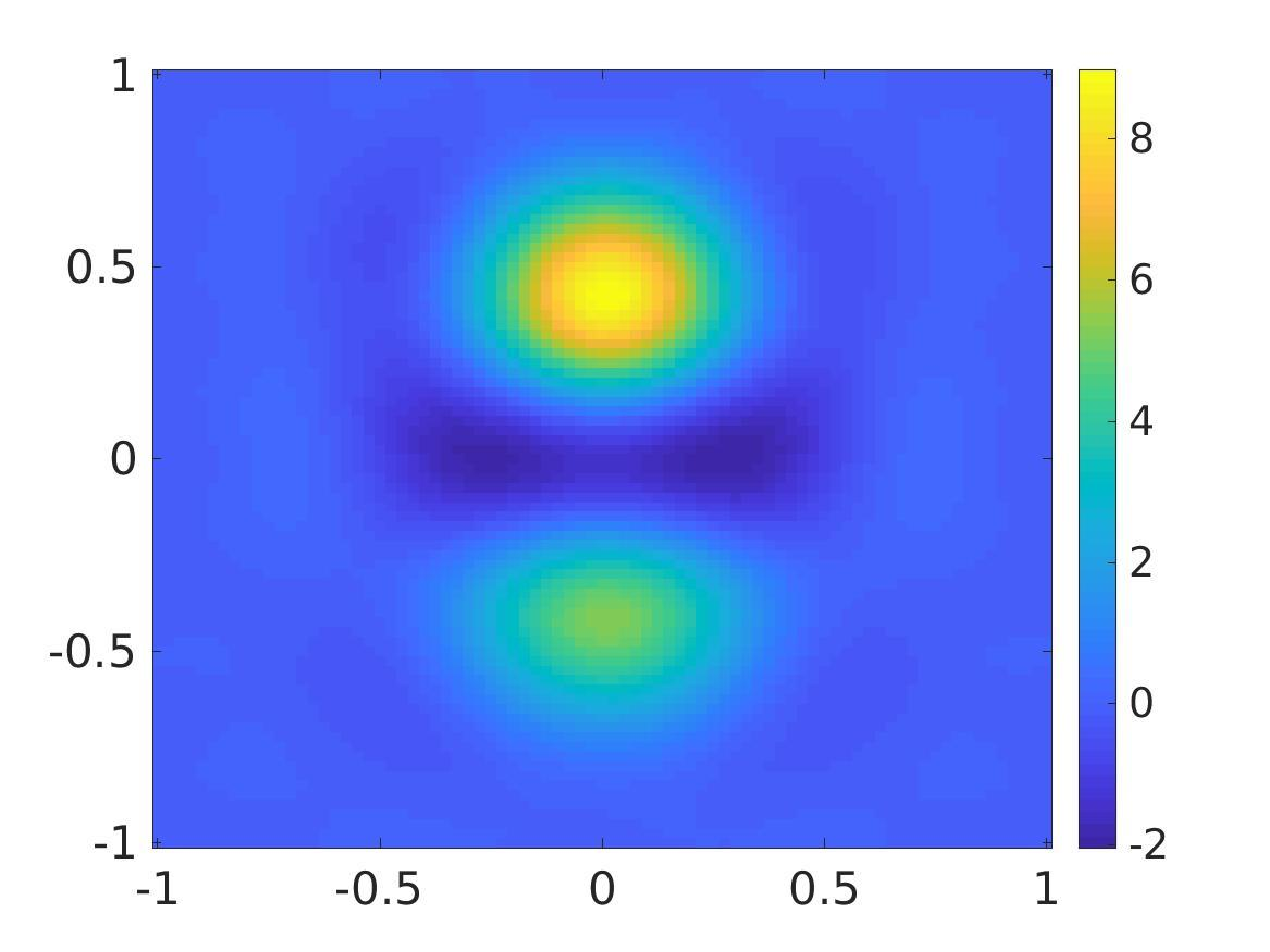}
	}		
	\subfloat[\label{3e}The error function ${E}$]{
		\includegraphics[width=.3\textwidth]{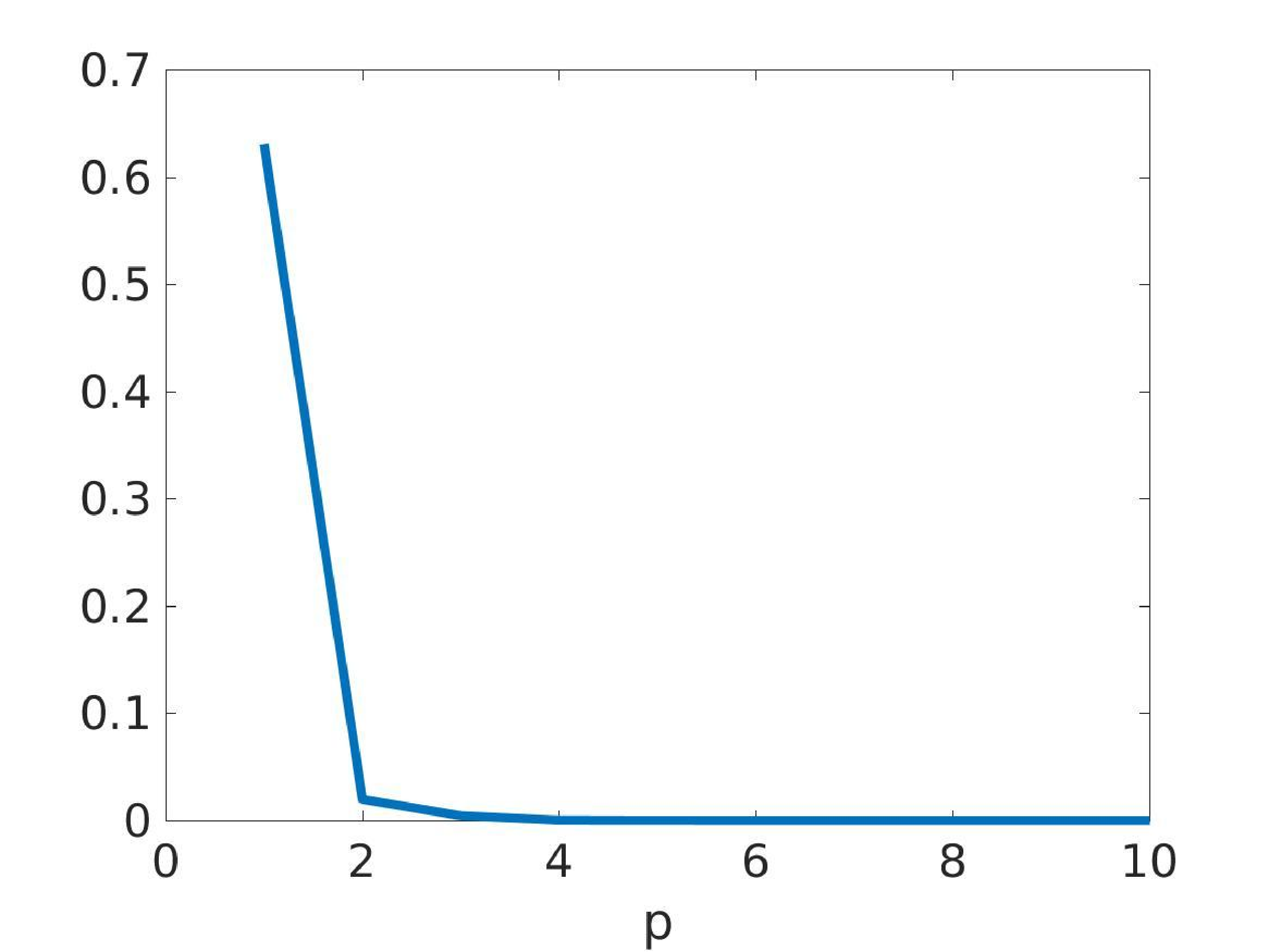}
	}

	\caption{\label{fig 3}Test 3. The true coefficient  and computed coefficient $c$. We already see both inclusions  in the graph of the first approximation $c^{(0)}$, computed by Step \ref{reg 1} of Algorithm \ref{alg}, see Figure \ref{3b}.  
	It is evident from the graph of the error function ${E}$, see Figure \ref{3e}, that the sequence $\{c^{(p)}\}_{{p \geq 1}}$ converges fast.
	}
\end{figure}

The numerical results for this case {are} displayed in Figure \ref{fig 3}. 
One can observe in Figure \ref{3b} that the circular shape and can be successfully detected at the first step.
The true maximal value of the function $c_{\rm true}$ at the lower inclusion is 8 and the reconstructed one is 8.90. 
 The true maximal value of the function $c_{\rm true}$ at the upper inclusion is 5 and the reconstructed one is 5.24. 
 Figure \ref{3e} shows the stability of our method.

\item {\it Test 4.} 
We test the case when the function $c_{\rm true}$ is allowed to be negative. 
In this case, the function $c_{\rm comp}$ is the letter $X$ with a half is positive and another half is negative.
The function $c_{\rm true}$ is given by
\[
	c_{\rm true}{(x, y)} = \left\{
		\begin{array}{ll}
			8 &|x| < 0.8, -0.8 < y \leq 0, |x \pm y| < 0.25,\\
			-8 &|x| < 0.8, 0< y < 0.8, |x \pm y| < 0.25,\\
			0 &\mbox{otherwise}.
		\end{array}
	\right.
\]

\begin{figure}[h!]
	\subfloat[The function $c_{\rm true}$]{
		\includegraphics[width=.3\textwidth]{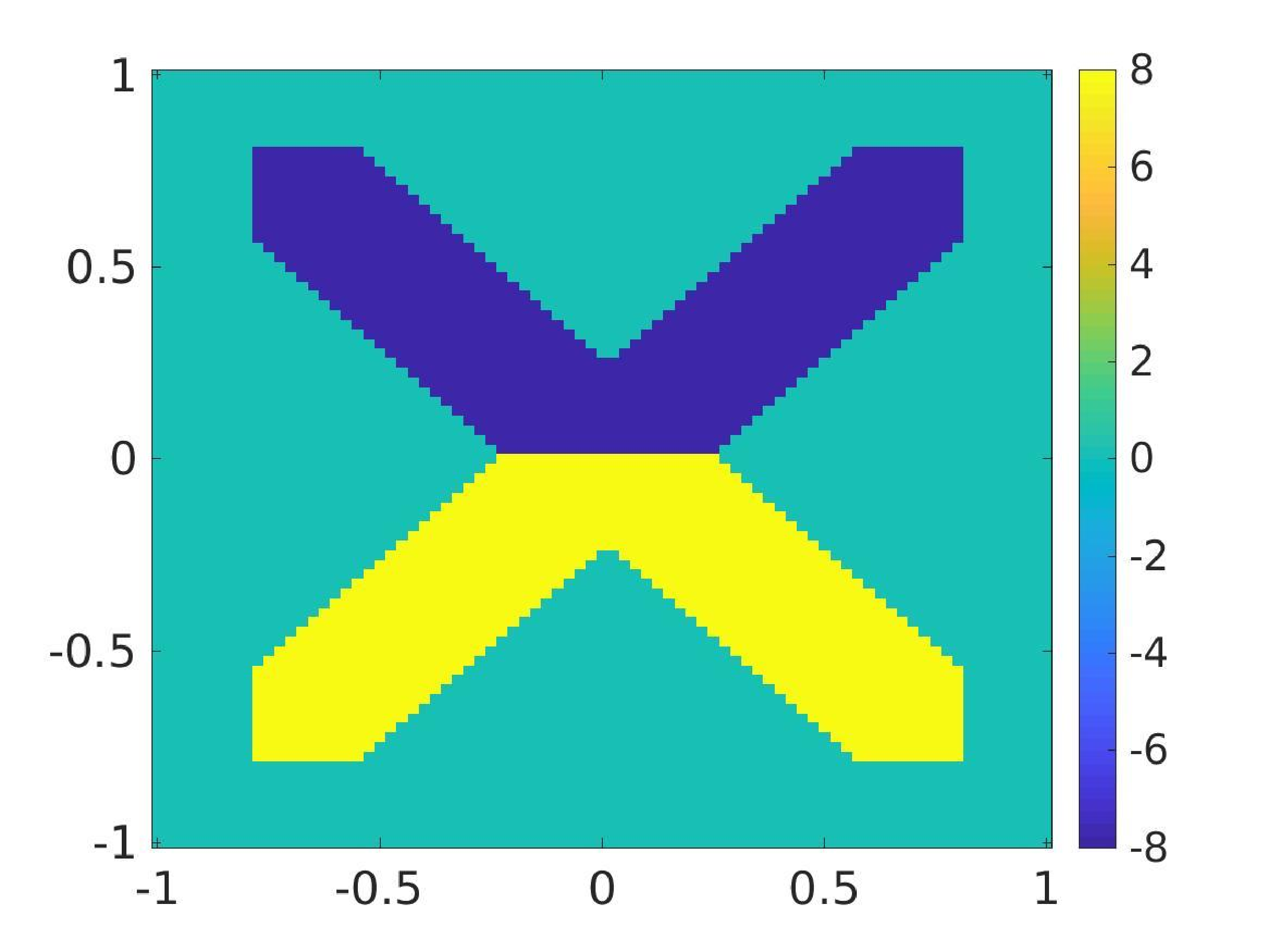}
	}  	
	\subfloat[\label{4b}The function $c^{(0)}$]{
		\includegraphics[width=.3\textwidth]{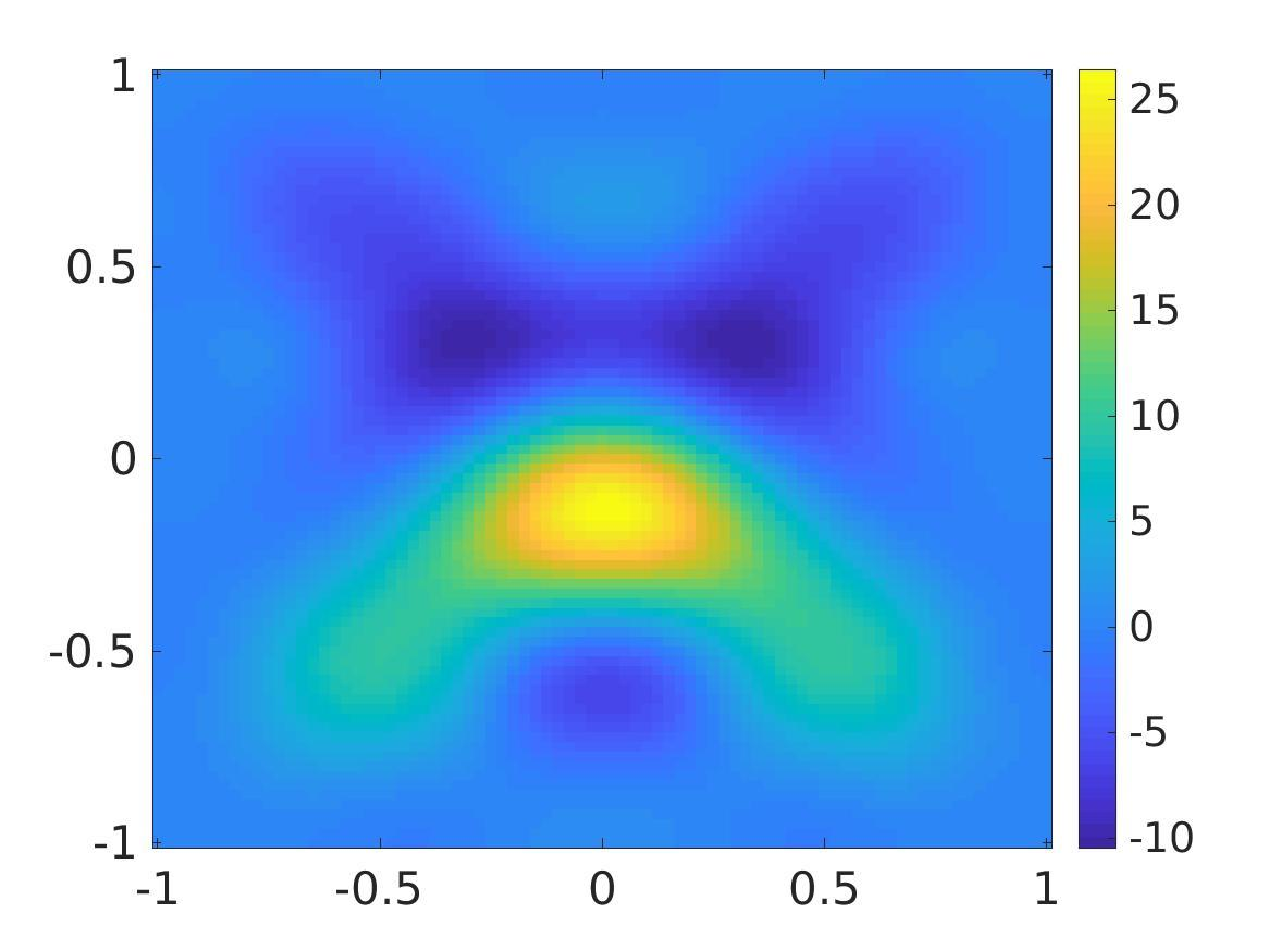}
	} 
	\subfloat[\label{4c}The function $c^{(1)}$]{
		\includegraphics[width=.3\textwidth]{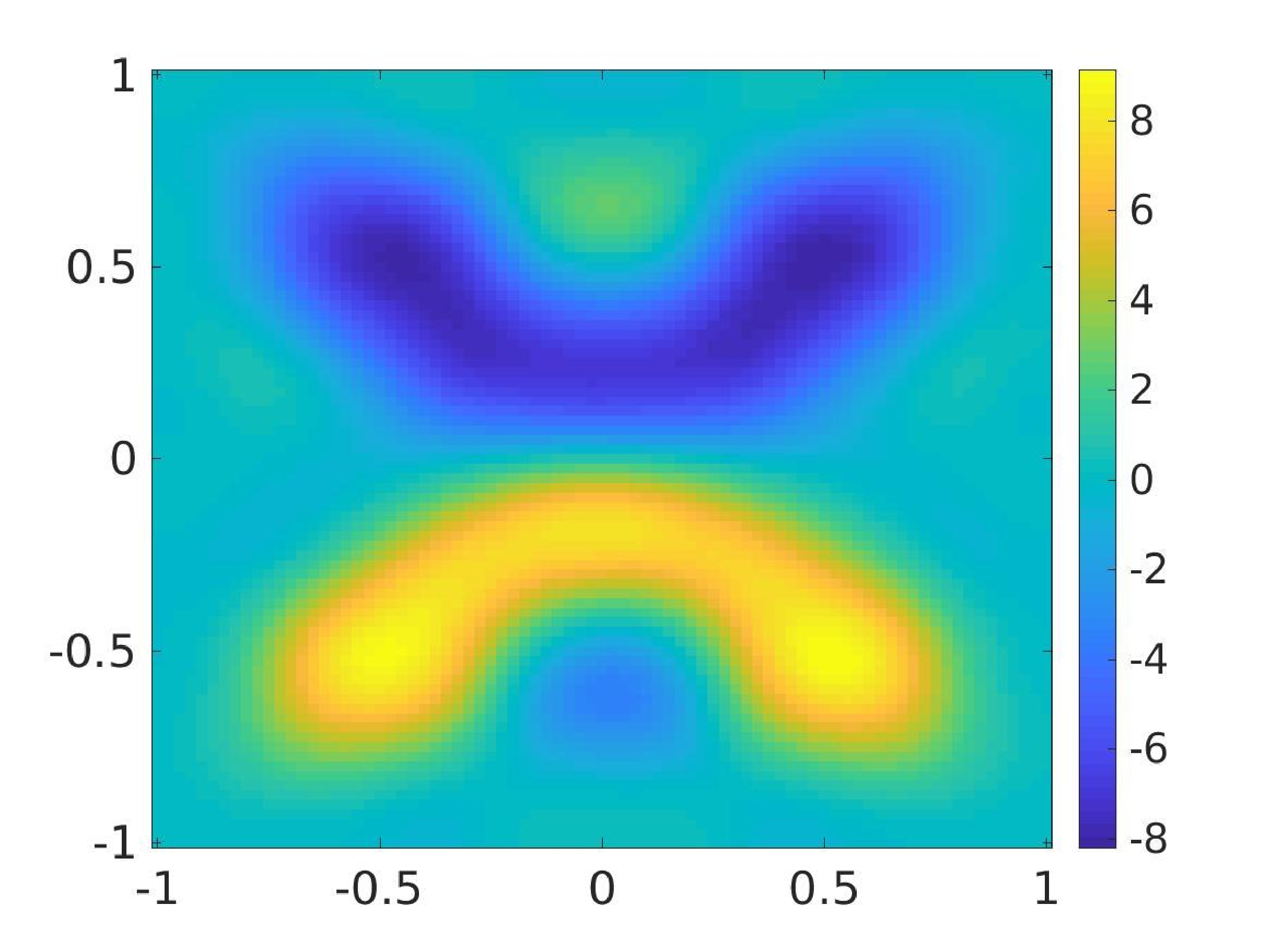}
		}	
	
	\subfloat[The function $c^{(2)}$]{
		\includegraphics[width=.3\textwidth]{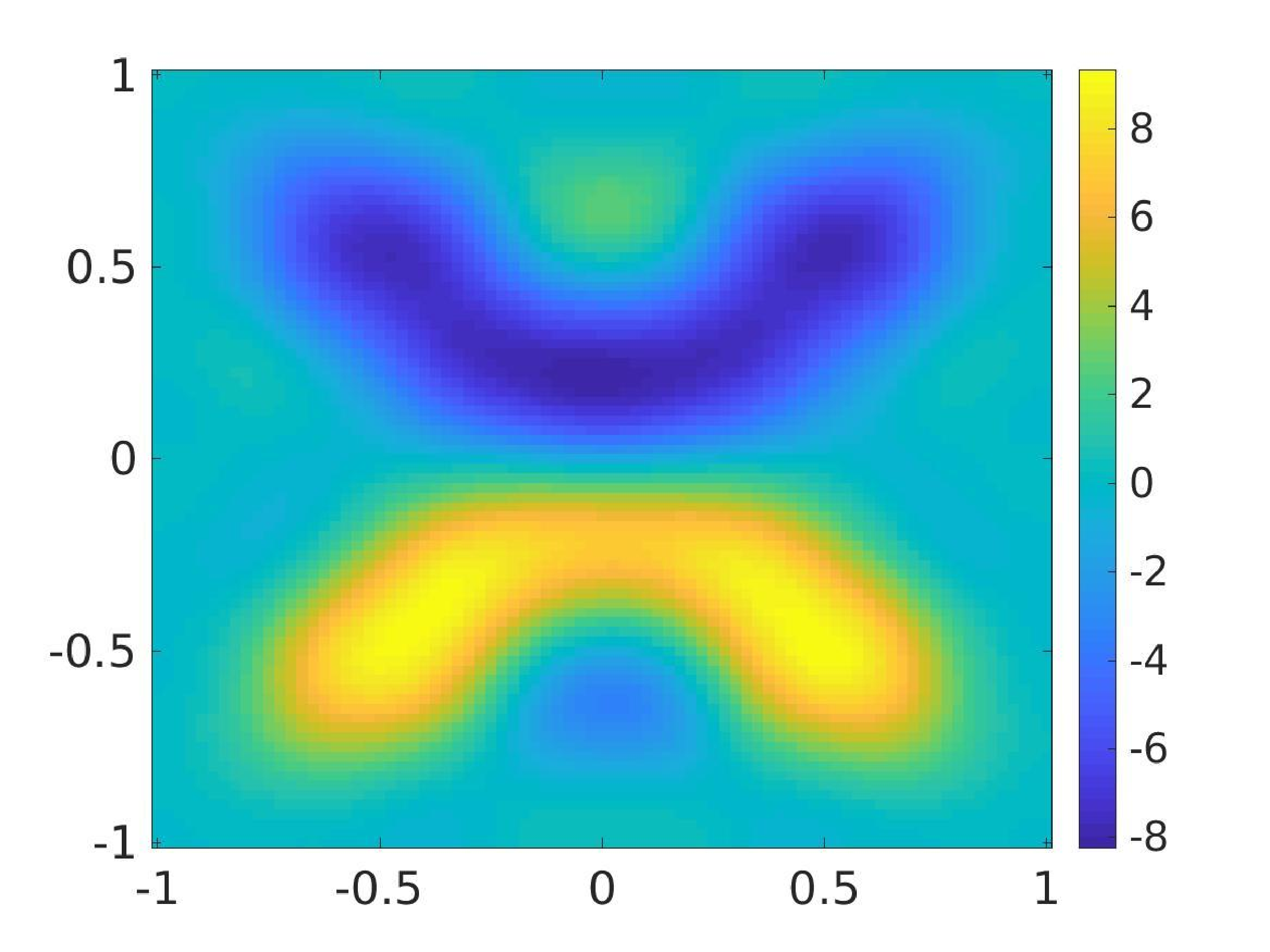}
	}	
	\subfloat[\label{4d}The function $c^{(10)}$]{
		\includegraphics[width=.3\textwidth]{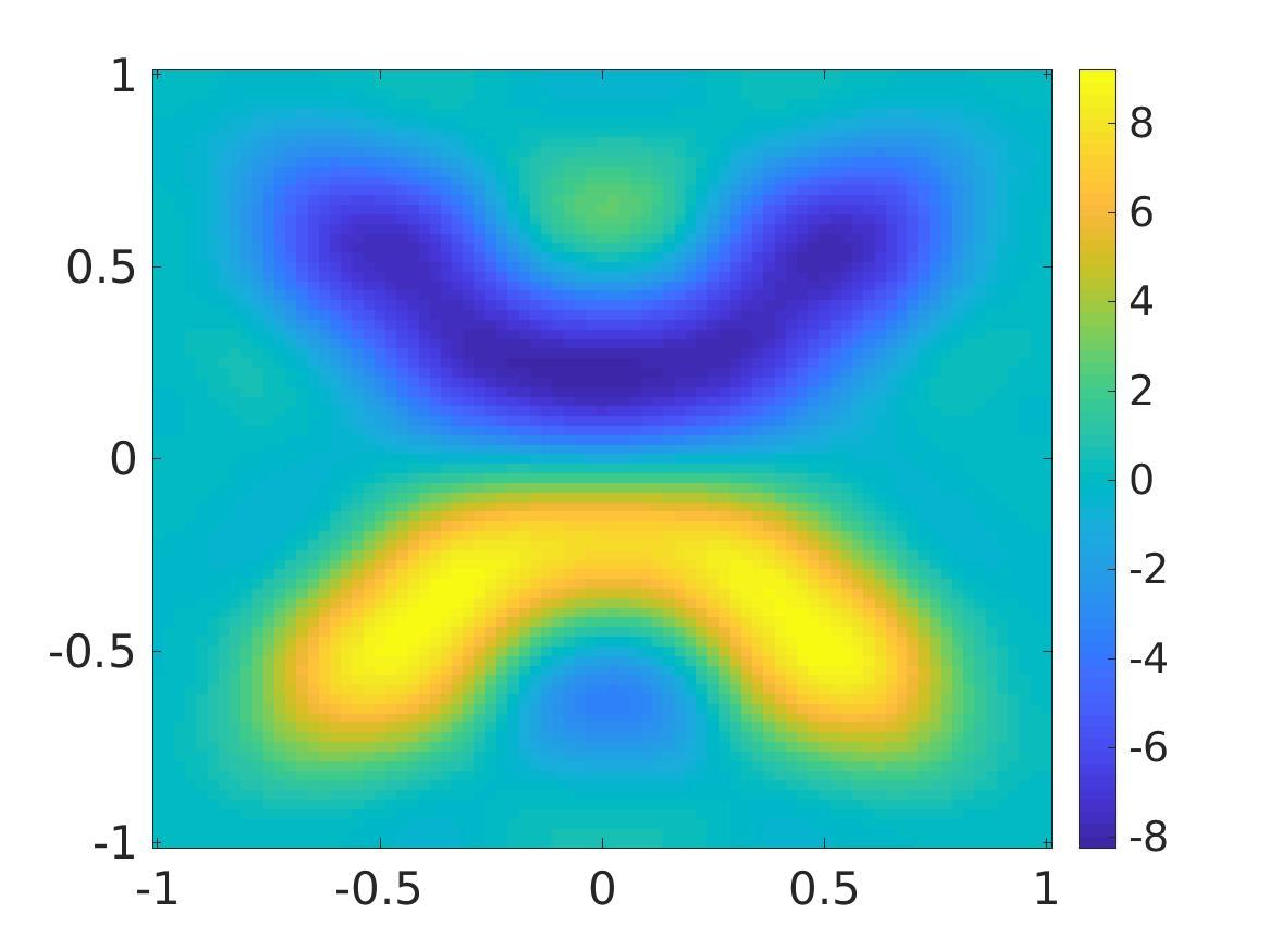}
	}		
	\subfloat[\label{4e}The error function ${E}$]{
		\includegraphics[width=.3\textwidth]{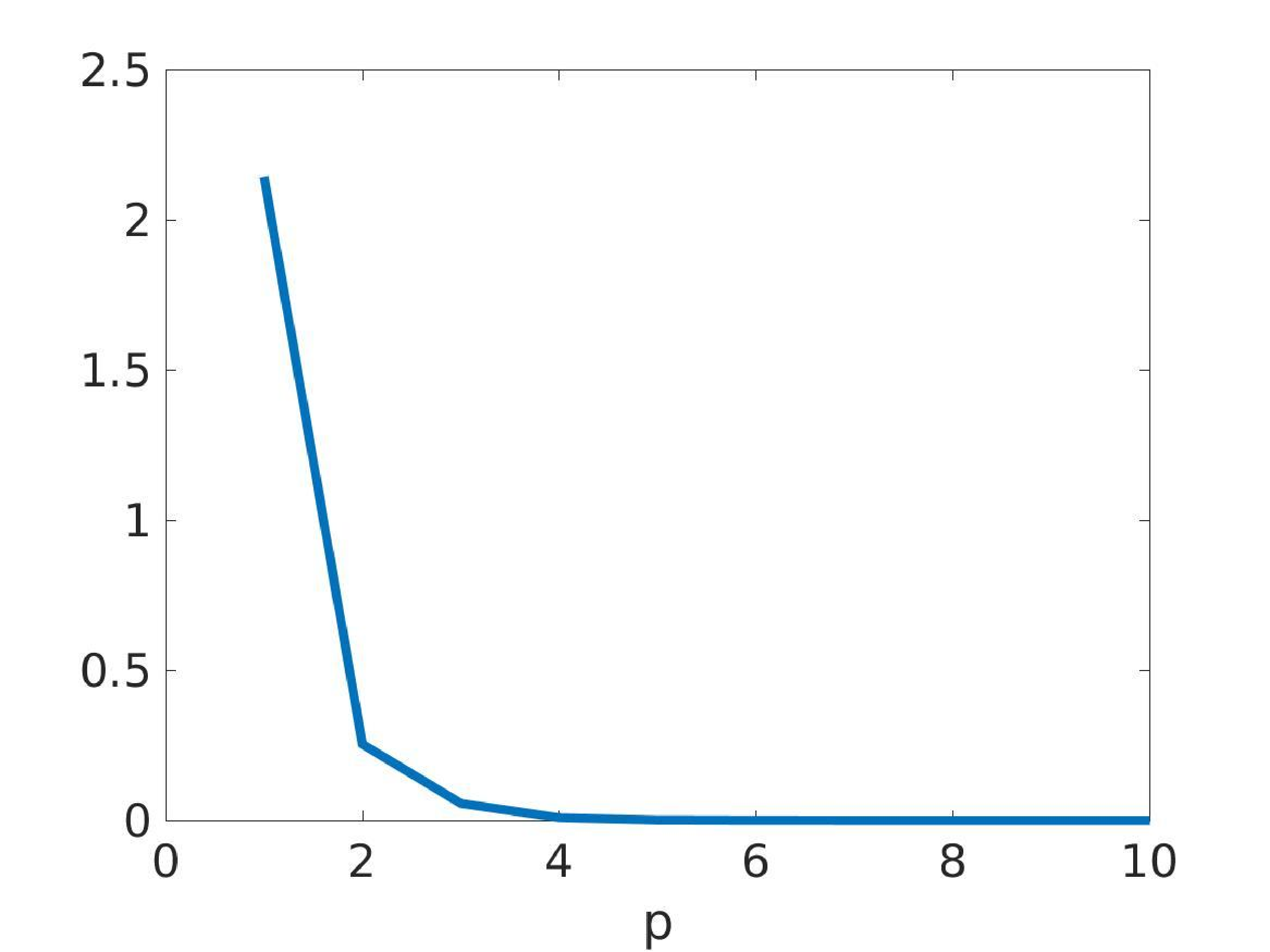}
	}

	\caption{\label{fig 4}Test 4. The true coefficient  and computed coefficient $c$. We already see the letter ``$X$"   in the graph of the first approximation $c^{(0)}$, computed by Step \ref{reg 1} of Algorithm \ref{alg}, see Figure \ref{4b}.  
	It is evident from the graph of the error function ${E}$, see Figure \ref{4e}, that the sequence $\{c^{(p)}\}_{{p \geq 1}}$ converges fast.
	}
\end{figure}

The numerical results for this case {are} displayed in Figure \ref{fig 4}. 
The reconstructed image of the letter $X$ is acceptable.
The true maximal positive value of the function $c_{\rm true}$  is 8 and the reconstructed one is 8.90. 
 The true minimal negative value of the function $c_{\rm true}$  is -8 and the reconstructed one is $-7.93$. 
Again, Figure \ref{4e} shows the stability of our method.

\end{enumerate}
{
\begin{remark}
	It is evident from Figures \ref{1c}--\ref{4c} that Algorithm \ref{alg} is robust in the sense that it provides good reconstructed coefficient $c_{\rm comp}$ after a few iterations without any requirement of an initial guess.
	It is remarkable mentioning that, in all tests above, our method provides good numerical results  without any advanced knowledge of the true coefficient $c_{\rm true}$.
	However, as seen in   \ref{1d}--\ref{4d}, there are some artifacts. 
	These artifacts might be caused  by cutting off the Fourier series of the function $v(\x, t)$ in \eqref{2.6}.
\end{remark}
}

\section{Concluding remarks} \label{sec concluding}

{In this paper, we introduced a new approach to numerically compute the creation or depletion coefficient of a general parabolic equation from lateral Cauchy data.
Although this problem is highly nonlinear, we successfully compute this coefficient without requiring a good initial guess.
In the first step, we derive an equation without the presence of the unknown coefficient. 
Then, we consider an approximation of the solution to this equation by truncating its Fourier series with respect to a special orthonormal basis of $L^2$.
By this, we obtain a system of nonlinear elliptic equations.
Numerically solving this system by an iterative procedure directly yields the desired coefficient.
Two of the important strengths of this method are that, unlike the optimal control method, we do not require a good initial guess for the true solution and that our algorithm converges fast. 
However, the main drawback in this paper is the analysis of the convergence of the sequence obtained by the algorithm.
A theorem that guarantees the efficiency of our method is not yet available.
That means, our method is verified only in the numerical point of view.
On the other hand,
we strongly believe that our technique can be applied to compute the diffusion coefficients of parabolic equations.
This serves as our near future work.}

\section*{Acknowledgments}
The author is grateful to Michael V. Klibanov for many fruitful discussions.

\end{document}